\documentclass[aos,MSNbibl,nameyear,seceqn,dvips]{arximspdf}
\usepackage{mathbh}
\usepackage{graphicx}

\doi{10.1214/12-AOS1040} 
\volume{40}
\issue{5}
\pubyear{2012}
\firstpage{2511}
\lastpage{2540}

\makeatletter

\newcommand{\rrvert}{\vert}
\newcommand{\llvert}{\vert}
\renewcommand{\ne}{\mathrm{NE}}

\newcommand{\pe}{\iota_n}
\newcommand{\implies}{\Longrightarrow}

\newcommand{\pin}{\pi_1,\ldots, \pi_n}
\newcommand{\Yij}{Y_{ij}}
\newcommand{\Xij}{X_{ij}}
\newcommand{\Zij}{Z_{ij}}
\newcommand{\M}{\mathrm{M}}
\newcommand{\W}{\mathrm{W}}
\newcommand{\Z}{\mathrm{Z}}
\newcommand{\C}{\mathrm{C}}
\newcommand{\I}{\mathrm{I}}
\newcommand{\A}{\mathrm{A}}
\newcommand{\R}{\mathrm{R}}
\newcommand{\D}{\mathrm{D}}
\newcommand{\Y}{\mathrm{Y}}
\newcommand{\olW}{\overline{\W}}
\newcommand{\eps}{\varepsilon}

\newcommand{\raw}{\rightarrow}

\newcommand{\E}{\mathrm{E}}
\newcommand{\pb}{{P}}

\newtheorem{theorem}{Theorem}[section]
\newtheorem{lemma}[theorem]{Lemma}

\newproclaim{remark}{Remark}[section]

\newproclaim{conditions}{Conditions}

\newcommand{\dm}{}

\newcommand{\hb}{{\hat{b}}}
\newcommand{\bW}{\bar{W}}
\newcommand{\bX}{\bar{X}}
\newcommand{\bY}{\bar{Y}}
\newcommand{\bbfY}{\bar{\mathbf{Y}}}
\newcommand{\bZ}{\bar{Z}}
\newcommand{\bbn}{\mathbb{N}}
\newcommand{\bbr}{\mathbb{R}}

\newcommand{\bbz}{\mathbb{Z}}
\newcommand{\F}{\mathcal{F}}
\newcommand{\cI}{\mathcal{I}}
\newcommand{\cX}{\mathcal{X}}

\newcommand{\al}{\alpha}
\newcommand{\be}{\beta}
\newcommand{\De}{\Delta}
\newcommand{\ep}{\varepsilon}
\newcommand{\Ga}{\Gamma}
\newcommand{\La}{\Lambda}
\newcommand{\ka}{\kappa}
\newcommand{\ph}{\phi}
\newcommand{\si}{\sigma}
\newcommand{\Up}{\Upsilon}

\newcommand{\hG}{\hat{G}}

\newcommand{\hal}{{\hat{\alpha}}}

\newcommand{\tX}{\tilde{X}}

\newcommand{\Var}{\operatorname{Var}}

\newcommand{\ind}{\mathbh{1}}

\newcommand{\dd}{d}

\newcommand{\nti}{n\rightarrow\infty}
\newcommand{\rai}{\rightarrow\infty}

\newcommand{\rl}{\mathbb{R}}
\newcommand{\eqref}[1]{(\ref{#1})}
\makeatother

\begin{document}
\begin{frontmatter}

\title{A penalized empirical likelihood method in high
dimensions\thanksref{T1}}
\runtitle{Penalized empirical likelihood}
\thankstext{T1}{Supported in part by NSF Grants DMS-07-07139
and DMS-10-07703.}

\begin{aug}
\author{\fnms{Soumendra N.}~\snm{Lahiri}\corref{}\ead[label=e1]{snlahiri@stat.tamu.edu}}
\and
\author{\fnms{Subhadeep}~\snm{Mukhopadhyay}\ead[label=e2]{deep@stat.tamu.edu}}
\runauthor{S.~N. Lahiri and S. Mukhopadhyay}
\affiliation{Texas A\&M University}
\address{Department of Statistics\\
Texas A\&M University\\
TAMU MS 3143\\
College Station, Texas 77843-3143\\
USA\\
\printead{e1}\\
\phantom{E-mail:\ }\printead*{e2}} 
\end{aug}

\received{\smonth{9} \syear{2011}}
\revised{\smonth{8} \syear{2012}}

%
\begin{abstract}
This paper formulates a penalized empirical likelihood
(PEL) method for inference on the population mean when the dimension
of the observations may grow faster than the sample size.
Asymptotic distributions of the PEL ratio statistic is derived
under different component-wise dependence structures of the
observations, namely, (i) non-Ergodic, (ii) long-range dependence
and (iii) short-range dependence. It follows that the limit
distribution of the proposed PEL ratio statistic can
vary widely depending on the correlation structure,
and it is typically different from the usual chi-squared limit
of the empirical likelihood ratio statistic in the fixed and
finite dimensional case. A~unified subsampling based calibration is
proposed, and its validity is established in all three
cases, (i)--(iii).
Finite sample properties of the method are investigated
through a simulation study.
\end{abstract}

%
\begin{keyword}[class=AMS]
\kwd[Primary ]{62G09}
\kwd[; secondary ]{62G20}
\kwd{62E20}
\end{keyword}

\begin{keyword}
\kwd{Asymptotic distribution}
\kwd{long-range dependence}
\kwd{Rosenblatt process}
\kwd{Wiener--It\^o integral}
\kwd{regularization}
\kwd{simultaneous tests}
\kwd{subsampling}
\end{keyword}

\end{frontmatter}

\section{Introduction}\label{sec1}
In a seminal paper, \citet{O88} introduced the empirical likelihood
(EL) method
for statistical inference on population parameters
in a nonparametric framework, and showed that it enjoyed
properties similar to the likelihood-based inference
methods in a more traditional parametric framework.
Following \citet{O88}, the EL method has been extended to various complex
inference problems; see, for example, Diccicio, Hall
and Romano (\citeyear{DHR91}), Hall and Chen (\citeyear{HC93}), \citet{qin94},
Owen (\citeyear{owen01}), \citet{Ber06}, Hjort, \mbox{McKeague} and Van Keilegom (\citeyear{H09}),
\citet{C09} and the references therein.
An extension of the EL method in the high-dimensional
context, where the dimension $p$ of the observations
increases with the sample size~$n$,
is given by Hjort, McKeague and Van~Keilegom (\citeyear{H09}).
Hjort, McKeague and Van~Keilegom (\citeyear{H09}) derives the limit distribution of
the EL ratio statistic based on $p$-dimensional estimating
equations when $p\raw\infty$ with $n$ at the rate $p= o(n^{1/3})$.
\citet{C09} improved upon the rate restriction in
Hjort, McKeague and Van~Keilegom (\citeyear{H09}) and established a nondegenerate
limit distribution of the EL ratio statistic,
allowing $p = o(n^{1/2})$ under suitable regularity conditions.

For applications to high-dimensional problems, such as
those involving gene expression data, one encounters
a $p$ that is typically much larger than the sample size $n$.
However, extension of the EL to such high-dimensional problems is
itself a daunting task because the (standard)
EL method is known to fail in such situations. An important
result of \citet{Tsao04} shows that the definition of \textit{the
EL for a $p$-dimensional population mean based on a sample size $n$
breaks down on a set of positive probability
whenever $p>n/2$}; further, this probability is asymptotically
nonnegligible. The main reason for this surprising
behavior of the EL is that for $p>n/2$,
the convex hull of $n$ random vectors in $\bbr^p$ is too small a
set to contain the true mean with high probability. As a
result, the standard EL approach cannot be applied to
the ``large $p$ small $n$'' problems with $p>n/2$. An alternative
formulation of the EL in such situations
(called the \textit{adjusted EL})
is given by \citet{CVA08}, which is
further refined and studied by \citet{EO09}.
The adjusted EL method adds additional pseudo-observations
[a single one in Chen, Variyath and Abraham (\citeyear{CVA08}) and two in \citet{EO09}]
so as to cover a hypothesized value of the mean parameter
within the convex hull of the augmented data set, thereby
making the adjusted EL well-defined. A second approach, due to
\citet{Bar07}, is to drop the convex hull constraint
in the formulation of the EL altogether and redefine the
likelihood of a hypothesized value of the parameter by
penalizing the unconstrained EL using
the Mahalanobis distance.
The penalized EL (PEL)
of \citet{Bar07} is well defined
for all values of $p\leq n$, as long as the sample covariance
matrix is nonsingular.
However, due to the use
of the inverse of the sample covariance matrix
in its formulation, the
PEL of \citet{Bar07} is also not well defined for
$p>n$. \citet{Bar07} establishes a chi-squared limit
of the PEL for the population mean in the case where the
dimension $p$ is fixed and finite for all $n$.
Other variants of the PEL where
\textit{a penalty function
is added to the standard EL}, in the spirit of the
penalized likelihood work of \citet{FL01}
and \citet{FP04},
are considered by \citet{Ot07}
and
\citet{TL10}. Both these papers consider the high dimensional
set up
and establish validity of their methods
still requiring $p$ to grow at most as
a fractional power of the sample size $n$.
\textit{In this paper, we introduce a modified version of the
PEL method of \citeauthor{Bar07} \textup{(\citeyear{Bar07})} that is computationally
simpler and that is applicable to a large class
of ``large $p$ small $n$'' problems, allowing $p$ to grow
faster than $n$}. This is an important step in generalizing
the EL in high dimensions
beyond the $p\leq n$ threshold where the standard EL
and its existing variants fail.

To briefly describe the proposed methodology and the
main results of the paper, suppose that
$X_1,\ldots, X_n$ are independent and identically distributed (i.i.d.)
$\bbr^p$-valued random vectors with mean $\mu\in\bbr^p$,
$1<p<\infty$.
Denote the $j$th component of a $p$-vector $x$ by $x_{j}$,
$j=1,\ldots,p$. The proposed PEL employs a multiplicative
penalty term to penalize
the likelihood of a hypothesized value $\mu$ of
the population mean as a quadratic function of
the distance between the sample mean and 
$\mu$. However, unlike Bartolucci's (\citeyear{Bar07}) method,
the use of the inverse sample covariance matrix is completely avoided,
as consistency of the sample covariance matrix in the high
dimensional case for all the dependence structures that we consider
in this paper is not guaranteed. The proposed PEL instead
uses a \textit{component-wise scaling} to bring up the
varying degrees of variability (variances) along different
components to a common level, and then it applies an overall
penalty on the sum of
the squared rescaled differences; see
\eqref{pel} in Section~\ref{sec2} below.
As a result, \textit{the proposed PEL is well-defined {for all}
values of $n,p\geq1$}. Further, this approach has the added
advantage that it does not require inversion of a high-dimensional
matrix, and therefore,
it is computationally much simpler.

For investigations into the theoretical properties of the
proposed PEL method,
we allow the components of $X_1$ to be dependent.
The range of dependence that we
consider covers the cases of:
\begin{longlist}[(iii)]
\item[(i)]  \textit{short-range dependence } (SRD), where
roughly speaking, the average of
the components of $X_1$ satisfies a central limit theorem (CLT)
under suitable
moment conditions; cf. \citet{IL09};
\item[(ii)]
\textit{long-range dependence} (LRD), where under suitable
regularity conditions, the average of the
components satisfies noncentral limit theorems
[Taqqu (\citeyear{T75,T77}), \citet{D79}];
\item[(iii)]
\textit{nonergodicity} (NE), where the dependence is so
strong that the average of the components even
fails to satisfy a (strong) law
of large numbers.
\end{longlist}
We refer to the LRD and SRD cases collectively
as the ergodic (E)-case, as the negative logarithm of
the PEL ratio
statistic $K_n$ (say) here satisfies a law
of large numbers without further centering and scaling,
for \textit{any} rate of growth of $p$; cf. Remark~4.2 below.
However, such degenerate limits laws are not always the
most useful in practice
as these only lead to conservative large sample inference procedures.
By using suitable centering and scaling, we are able to further refine
these results and
establish convergence in distribution to nondegenerate limits.
Specifically, we show that under SRD,
$K_n$ with centering at 1
[for $c_*=1$ in condition~(C.2)(ii) below]
and scaling by square-root of the dimension $p$ of the observations
converges to a Normal limit, very much like the results of
Hjort, McKeague and Van~Keilegom (\citeyear{H09}) and \citet{C09}, but allowing
a much faster rate of growth of $p$ and allowing a more general
dependence framework.
In the long range dependent (also abbreviated as LRD) case,
%
$K_n$ with a suitable normalization can have both
Normal and non-Normal limits. For the Normal limit, the
centering and the scaling sequences are the same as those
used in the SRD case, except at the boundary layer
of dependence where the Normal limit switches over to
the non-Normal limit. For the non-Normal limit under LRD,
the centering term is the same as that in the SRD case, but the scaling
depends on the rate of decay of the auto-correlation coefficient
of the components of $X_1$ (up to a possibly unknown permutation).
Finally, in comparison to the E-case,
$K_n$
in the NE-case is shown to converge
in distribution to a stochastic integral, and it
does NOT require any further centering and scaling.

%
\begin{figure}

\includegraphics{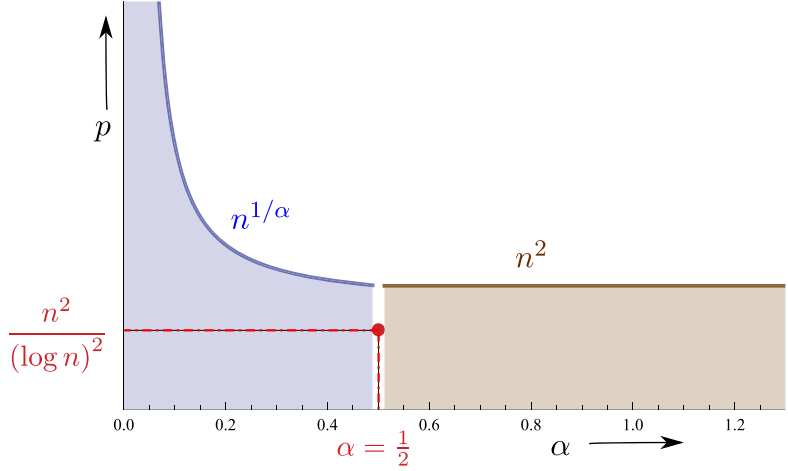}

\caption{Envelope for the growth rates of $p$ for a nondegenerate
limit law of the PEL ratio statistic. The smaller the $\al$, the
stronger is the dependence.}\label{fig01}
\end{figure}

%

%

The growth rate of $p$, for which a \textit{nondegenerate limit
law} holds for a suitably transformed $K_n$,
primarily depends on the strength of dependence among
the components of the observations; cf. Figure~\ref{fig01}.
In the NE-case, $p$ can grow
\textit{arbitrarily fast} (e.g., polynomial, exponential,
super-exponential,
etc.)
as a function of $n$. In the E-case, although a degenerate limit
law holds for an arbitrary growth rate of~$p$, for a nondegenerate limit,
$p$ must admit a suitable upper bound. In particular,
for the Normal limit in the E-case (excluding the boundary case),
the growth rate of $p$ as a function of $n$ is
$p=o(n^2)$. For the non-Normal limit under the E-case,
$p= o(n^{1/\al})$ for $0<\al<1/2$
where, roughly speaking, $\al$ denotes the exponent
of the rate of decay of the autocorrelation among the
components of $X_1$, up to a permutation; cf. condition (C.4)$_\al$,
Section~\ref{sec3}.
The boundary case is given by $\al=1/2$, where the
growth rate is slightly smaller and is given by
$p = o(n/[\log n]^2)$.
From Figure~\ref{fig01}, it follows that
the stronger the dependence among the
components of the observations, the higher is the
allowable growth rate of $p$ as a function of $n$
for a nondegenerate limit. The limiting case $\al\raw0+$ is the
NE-case. Here the nondegenerate limit for the
the negative logarithm of the PEL ratio statistic
holds for an arbitrary growth rate
of $p$ as a function
of $n$.

It is worth pointing out that in most cases, the
limit distribution of the PEL ratio statistic is \textit{not}
distribution free in the sense that the asymptotic
approximation to the distribution of the
PEL ratio statistic requires the knowledge of one or more
unknown population parameters. As a result,
the limit laws are not directly usable in practice.
To address this issue, we propose a calibration procedure based on the
subsampling method. We show that under mild conditions,
the subsampling based calibration method
is consistent
under \textit{all three}
types of dependence structures.

The key step in the proofs is to derive a quadratic
asymptotic approximation to $K_n$
%
under all three cases of dependence. This is
presented in Lemma~\ref{lem-2} for the NE-case and in the proof of
Theorem~\ref{th3.2} for the E-case.
The derivation of the limit law in the NE-case
uses some weak convergence and operator convergence results
on Hilbert spaces.
On the other hand, in the E-case,
refined approximations to $K_n$
are required to go beyond their degenerate limits.
See Section~\ref{sec6} for more details.


The rest of the paper is organized as follows. In Section~\ref{sec2}, we
describe the
PEL methodology. In Section~\ref{sec3}, we introduce the asymptotic framework
and establish the limit distributions of the logarithm of
the PEL ratio
statistic under the three dependence scenarios. In Section~\ref{sec4}, we
describe the subsampling method and prove its validity for
all three cases. We report the results
from a moderately large simulation study in Section~\ref{sec5}.
Proofs of the results are given in Section~\ref{sec6}.

\section{Formulation of the PEL}\label{sec2}
Let $X_1,\ldots,X_n$ be i.i.d. random vectors with
mean $\mu\in\bbr^p$. Let $X_{ij}$ denote the $j$th component of
$X_i$, $1\leq i\leq n$, $1\leq j\leq p$.
Also, let $A'$ denote the transpose of
a matrix $A$. We define the \textit{penalized empirical likelihood (PEL)}
of a plausible value $\mu=(\mu_1,\ldots,\mu_p)'$
of the population mean as
%
\begin{equation}
L_n(\mu) =\sup_{(\pi_1,\ldots,\pi_n)'\in\Pi_n} \Biggl\{ \Biggl(\prod
_{i=1}^n \pi_i \Biggr) \exp \Biggl( -
\lambda\sum_{j=1}^p \delta_j
\Biggl[ \sum_{i=1}^n\pi_i
(X_{ij} - \mu_j) \Biggr]^2 \Biggr) 
\Biggr\}, \label{pel}\hspace*{-35pt}
\end{equation}
where $\Pi_n=\{(\pi_1,\ldots,\pi_n)'\in[0,1]^n$ :
$
\sum_{i=1}^n \pi_i =1\}$, $\delta_j$'s are component specific
weights (which may be random), and
$\lambda=\lambda_n\in[0,\infty)$ is an overall penalty factor.
Here we use
%
\begin{equation}
\delta_j\equiv\delta_{nj} = s_{nj}^{-2}
\ind(s_{nj}\neq0), \label{de-j}
\end{equation}
where $s_{nj}^2 = n^{-1}\sum_{i=1}^n(X_{ij} - \bX_{nj})^2$
is the sample variance of the $j$th components of
$X_1,\ldots, X_n$, $\bX_{nj}= n^{-1}\sum_{i=1}^nX_{ij}$
and where $\ind(\cdot)$ denotes the indicator function.
This choice of the component-wise
scaling allows us to adjust for the heteroscedasticity
along different co-ordinates of $X_1$ and therefore, the
overall penalization parameter $\lambda$
gives \textit{comparable} weights to \textit{all}
components. In addition, the choice of the penalty
function makes the proposed PEL
\textit{invariant} with respect to component-wise scaling,
which is an inherently desirable property, particularly
while dealing with high-dimensional variables, where
the assumption of homoscedasticity among a large number
of components is unrealistic.
The maximizer of the product $\prod_{i=1}^n \pi_i$ in \eqref{pel}
without the penalty term is given by $\pi_i=1/n$ for
$i=1,\ldots, n$. Hence, the \textit{PEL ratio statistic}
at a plausible value $\mu$ of the mean vector is
defined as
\[
R_n(\mu) = n^n L_n(\mu).
\]

We now compare our formulation with the PEL of
\citet{Bar07}, which is defined as
%
\begin{equation}\qquad
L^{\mathrm{B}}_n(\mu) = \sup_{(\pi_1,\ldots,\pi_n)'\in\Pi_n} \Biggl\{ \Biggl(\prod
_{i=1}^n \pi_i \Biggr) \exp
\biggl( - \frac{n (v-\mu)'V_n^{\dagger-1}(v-\mu)}{2h^2} \biggr) \Biggr\}, \label{b-pel}
\end{equation}
where $v=\sum_{i=1}^n\pi_i X_i$, $V_n^{\dagger}=
n^{-1}\sum_{i=1}^n(X_i -\bX_n)(\bX_i - \bX_n)'$ is the sample
covariance matrix, and $h$ is a penalty parameter.
Note that for
a large $p\in(1,n]$, the sample variance matrix $V_n$
is ill-conditioned, if the smallest eigen-value of the
underlying covariance matrix $\Sigma$ (say) of $X_1$ is not bounded
away from zero, and it is always singular
when $p>n$, requiring further modifications to make
\eqref{b-pel} well defined. Under either of the two
scenarios, the PEL based on \eqref{b-pel}
can be computationally demanding and unstable.
In comparison, the component-wise scaling in \eqref{pel}
only involves $1$-dimensional operations
which is computationally much simpler and feasible
even for a large $p$. A second limitation of \eqref{b-pel}
is the lack of attractive theoretical properties of
$V_n^{-1}$ (or of its variants) in high dimensions.
Indeed, consistency of the sample
covariance matrix (and its banded or tapered versions)
in the high-dimensional setting is questionable in presence
of strong correlations among the components of $X_1$
that we consider here. Existing work
on consistency of the sample covariance matrix
is known only under suitable conditions of sparsity or
weak dependence; cf. Bickel and Levina (\citeyear{bickel08}), \citet{kar08},
Cai, Zhang and Zhou (\citeyear{cai09}).
Our formulation also avoids this problem altogether by
using component-wise scaling.


For the sake of completeness, we also briefly describe
the penalized EL approach of \citet{Ot07} and \citet{TL10}, specialized to the case of the mean parameter
$\mu$ for simplicity of exposition. Let
$L_n^{\mathrm{ST}}(\mu)=
\sup\{\prod_{i=1}^n \pi_i
 \dvtx(\pi_1, \ldots,\pi_n)\in\Pi_n, \sum_{i=1}^n\pi_i (X_i-\mu
) =0
\}
$ denote the standard EL for $\mu$. Also, let $p_\lambda(\cdot)$ be a
penalty function, such as the smoothly clipped absolute deviation
(SCAD) penalty function of \citet{FL01}. Then, the penalized EL\vadjust{\goodbreak}
considered by \citet{Ot07} and \citet{TL10} is of the form
%
\begin{equation}
L_n^{\mathrm{OTL}}(\mu) = L_n^{\mathrm{ST}}(\mu)\exp
\Biggl(-n \sum_{j=1}^{p}p_\lambda(
\mu_j) \Biggr), \label{O-TL}
\end{equation}
where $\mu=(\mu_1,\ldots,\mu_p)'$. For the case of a
more general parameter $\theta$ defined through a
set of estimating equations,
the formulation of \citet{Ot07} and \citet{TL10} replaces $L_n^{\mathrm{ST}}(\mu)$ in \eqref{O-TL}
by the corresponding version of the standard EL for $\theta$; cf. \citet{qin94}.
As a result, irrespective of the target parameter, since \eqref{O-TL}
is directly based on the standard EL, this formulation of the
penalized EL also suffers from the same limitations as the standard EL.
In particular, this approach also fails in high dimensions
whenever $p>n/2$.

In the next section, we investigate theoretical properties of the
proposed PEL method \eqref{pel} under the dependence
structures described in Section~\ref{sec1}.

\section{Limit distributions}\label{sec3}
\subsection{General framework}\label{sec3.1}
We establish the limit distribution theory for the PEL
ratio statistic in a triangular array
set up, with $n$ denoting the variable driving the
asymptotics. Thus, the vectors $X_1,\ldots, X_n$
depend on $n$ as are their distributions
and the dimension $p$. However, we often
suppress the dependence on $n$
for simplicity of notation.
The limit distribution of the PEL
ratio statistic depends on the degree of dependence
among the components of $X_1,\ldots, X_n$.
As stated in Section~\ref{sec1},
we can broadly
classify the dependence structure into two categories: (i)~Non-Ergodic (NE) and (ii) Ergodic (E). In the NE-case, the
dependence among the components of each $X_i$ is so strong
[cf. condition (C.3) below] that
even the law of large numbers fails. In this case, we show that
under appropriate conditions, the PEL ratio
statistic has a nondegenerate limit. In contrast, in the E-case,
the corresponding limit is degenerate, and further centering and scaling
are needed for nondegenerate limit laws
which, in turn, depend the type of dependence (SRD or LRD).
We begin with the
NE-case.

\subsection{Limit distribution for the nonergodic case}\label{sec3.2}
We need to introduce some notation at this stage.
Let $\rho_n(j,l)$ denote the correlation between $X_{1j}$ and
$X_{1l}$, $1\leq j,l\leq p$. Let $\si_{nj}^2 = \operatorname{Var}(X_{1j})$
and $D^X_{nj} = \{x\in\bbr\dvtx P(X_{1j} = x)>0\}$, $1\leq j\leq p$,
$n\geq1$.
Write $C$ to denote a generic
constant in $(0,\infty)$.
Also, for any two sequences $\{a_n\}_{n\geq1}$ and
$\{b_n\}_{n\geq1}\in(0,\infty)$, write $a_n\sim b_n$
if $a_n/b_n\raw1$ as $\nti$. Let $L^2[0,1]$ denote
the set of all square integrable
functions on $[0,1]$ (with respect to the Lebesgue measure
on $[0,1]$), equipped with the inner product $\langle f,g \rangle
= \int_0^1 fg $, $f,g\in L^2[0,1]$. Let
$\{\phi_k \dvtx k\in\bbn\}$ denote
a complete orthonormal basis of $L^2[0,1]$, where $\bbn=\{1,2,\ldots\}$
denotes the set of all positive integers.
For any (bounded) jointly measurable function $h\dvtx[0,1]^2 \raw\bbr$,
define the operator $\Up_h$ on $L^2[0,1]$ by
$\Up_h f = \int_{[0,1]} h(\cdot;u)f(u)\,du$ for $f\in L^2[0,1]$.
For $x,y\in\bbr$, let $x\wedge y = \min\{x,y\}$.

We shall make use of the following
conditions for deriving the limit distribution of
$- \log R_n(\mu)$. The values of the integers $r$ and $s$
below will be specified later in the statements of the
theorems.

\begin{conditions*}
\begin{enumerate}[(C.3)]
\item[(C.1)] \hspace*{3pt}(i)  $\max\{ E[\si_{nj}^2\delta_{j}]^s\dvtx1\leq
j\leq
p\} = O(1)$
for a given $s \in\mathbb{N}$.
\begin{enumerate}[(ii)]
\item[(ii)]  $\limsup_{\nti}  \max\{P(X_{nj}=x)\dvtx x\in D^X_{nj}, 1\leq
j\leq
p\}
<1$.
\end{enumerate}
%
\item[(C.2)] \hspace*{3pt}(i)  For a given $r\in\bbn$,
$
\max\{E|X_{1j}|^{r}\dvtx1\leq j\leq p\} < C.$
\begin{enumerate}[(ii)]
\item[(ii)]  $\lambda_n = c_* n/p$ for some $c_*\in(0,\infty)$.
\end{enumerate}
\item[(C.3)] There exists a correlation function $\rho_0(\cdot,\cdot)$
of a mean-square
continuous 
process on $[0,1]$, and
for each $n\geq1$, there exists a permutation
$\pe$ of $\{1,\ldots, p_n\}$
such that
$\rho_n(j,l) =\rho_0 (\frac{\pe(j)}{p}, \frac{\pe(l)}{p}  )$.
Further, with $c_*$ as in (C.2),
$\rho_0(\cdot,\cdot)$ satisfies the following:
\begin{enumerate}[(iii)]
\item[(i)]  $4c_*^2 \int_0^1\int_0^1 \rho_0^2(u,v)\,du\,dv <1$;
\item[(ii)]
$
\sup\{ |\rho_0(u+h_1, v+h_2) - \rho_0(u, v)| \dvtx|h_1|\leq\delta,
|h_2|\leq\delta\}\leq
g(\delta) H(u,v)$
for all $u,v, u+h_1, v+h_2\in[0,1]$ for some function $g(\cdot)$ satisfying
$g(\delta)\raw0$ as $\delta\downarrow0$ and for some function
$H(\cdot, \cdot)$
satisfying
$
\sum_{k\geq1}  \langle|\phi_k|, |\Up_H \phi_k|  \rangle
<\infty;
$
\item[(iii)] $\sum_{k\geq1}  \langle\phi_k, \Up_0 \phi_k
\rangle$
converges, where $\Up_0=\Up_h$ with $h= \rho_0$.
\end{enumerate}
\end{enumerate}
\end{conditions*}

We now briefly comment on the conditions. Condition (C.1)(i)
is a moment condition on the scaled component-wise weights
$\delta_j$'s and requires finiteness of the $s$th negative
moment of the sample variance $s_{nj}^2$'s (scaled by
the respective expected values $\si_{nj}^2$'s).
This condition holds in the case of
Gaussian $X_{ij}$'s whenever the the sample size
$n >2s$. Condition (C.1)(ii) is a mild condition---it
says that none of the $X_{ij}$'s
take a single value with probability approaching one.
This condition trivially holds if the components
of $X_1$ are continuous and also in the discrete case, if the
supports of $X_{1j}$'s contain at least two values with
asymptotically nonvanishing
probabilities.
Condition (C.2)(i) is a moment condition that will be used
with different values of $r$ in the main theorems of
this section, while (C.2)(ii) specifies the growth
rate of the penalty parameter for a nondegenerate limit
of $-\log R_n(\mu_0)$.
However, unlike the standard usage of the penalty
parameter in the context of variable selection, where
different choices of the parameter lead to different sets
of variables being chosen, here the key role of the
penalty parameter is to stabilize the contribution
from the sum of component-wise squared differences
to the overall ``likelihood'' in \eqref{pel}. Finally, consider
condition (C.3) 
that specifies the nonergodic structure of the $X_i$'s.
Note that, up to a
(possibly unknown) permutation of the
co-ordinates, the components of $X_1$ are essentially correlated
as strongly as the variables $W(t)$'s coming from
a constant mean, mean-square continuous process $\{W(t) \dvtx t\in
[0,1]\}$ (say)
with covariance function $\rho_0(\cdot)$.
In this case, the dependence among the variables $W(i/p)$ and
$W(j/p)$, $1\leq i<j\leq p$ is so strong that the average
$p^{-1}\sum_{j=1}^p W(j/p)$ may not converge to a constant
as $p\raw\infty$, as one would expect from the well-known
ergodic theorems.

Under conditions (C.1)--(C.3),
the limit distribution of the log-PEL
ratio statistic is given by a stochastic
integral, as shown by the following result.

\begin{theorem}\label{th3.1}
Let conditions \textup{(C.1)}, \textup{(C.2)} and \textup{(C.3)} hold with $s=4$ and $r=8$,
let $\mu_0$ denote the true value of $\mu$
and let $p\raw\infty$ as $\nti$.
Then
%
\begin{equation}
- \log R_n(\mu_0) \raw^d c_*\int
_0^1 \int_0^1
\La_0(u,v) Z(u)Z(v) \,du\,dv, 
\label{v-str-dep}
\end{equation}
where $Z(\cdot)$ is a zero mean { Gaussian}
process on $[0,1]$ with
covariance function $\rho_0(\cdot)$ and where the function
$\La_0(\cdot,\cdot)$ is defined as
%
\begin{equation}
\La_0(u,v) = \sum_{k=0}^\infty(-2)^kc_*^k
\rho_0^{*(k)} (u,v), \qquad  0<u,v<1, \label{La-def}
\end{equation}
with
$\rho_0^{*(0)}(u,v) = 1$,   $\rho_0^{*(1)}(u,v) = \rho_0 (u,v)$
and for $k\geq1$,
\[
\rho_0^{*(k+1)}(u,v) = \int_0^1
\cdots\int_{0}^1 \Biggl\{\prod
_{j=1}^{k-1} \rho_0(u_j,
u_{j+1}) \Biggr\} \rho_0(u_1, u)
\rho_0(u_k, v) \,du_1\cdots
\,du_k.
\]
%
\end{theorem}

For a general definition of a stochastic integral of the form
\eqref{v-str-dep}, see Cramer and Leadbetter (\citeyear{cra67}). Note that
under condition (C.3), by repeated application of the Cauchy--Schwarz
inequality, for $k\geq1$,
\[
\sup_{u,v \in[0,1]} \bigl| \rho_0^{*(k+1)} (u,v)\bigr | \le(1-
\delta)^{k-1} [2c_*]^{-(k+1)}\qquad \mbox{for some }  \delta\in(0,1),
\]
and hence, 
the limiting stochastic integral is well defined.

Theorem~\ref{th3.1} shows that under a suitable choice of the penalty
parameter, namely, $\lambda=c_* n/p$, the negative log PEL ratio
statistic has a nondegenerate limit distribution. Note that,
unlike the standard version of the EL, we do not use the multiple $2$
before $- \log R_n(\mu_0)$. This is a direct artifact of
the additional penalty term that we use in the formulation
of PEL. Also, unlike most high-dimensional problems where the validity
of a large sample inference procedure
breaks down beyond a certain (often exponential)
rate of growth of $p$,
the PEL and the associated limit distribution of $- \log R_n(\mu_0)$
in the NE-case remains valid for \textit{arbitrary} rate of growth of $p$
as a function of the sample size. Thus, in the NE-case, it is possible
to carry out simultaneous hypothesis testing for a very large number of
parameters even with a moderately large sample.


\subsection{Limit distribution for the ergodic case}\label{sec3.3}
In the E-case, we shall make use of the following conditions,
for
$\al\in(0,\infty)$:

\begin{longlist}[(C.4)$_\al$]
\item[(C.4)$_\al$]
There exists a covariance function $\rho_\al(\cdot)$
on $\bbz\equiv\{0,\pm1,\pm2,\ldots\},$
and for each $n\geq1$, there exists a (possibly unknown)
permutation $\pe$ of $\{1,\ldots,
p_n\}$ such that $\rho_{\al}(k)\sim C |k|^{-\al}$ as $k\rai$
and
\[
\sup_{1\leq j,l\leq p} \biggl|\frac{\rho_n(j,l)}{\rho_\al (\pe(j) - \pe(l) )} -1 \biggr| = o(1)\qquad\mbox{as }
\nti. 
\]
\item[(C.5)$_\al$]
There exists a constant $C>0$ such that
\[
\check{\varrho}_n(k) \leq Ck^{-\al}\qquad \mbox{for all } k
\geq1, n>C,
\]
where $\check{\varrho}_n(\cdot)$ denotes the $\varrho$-mixing
coefficient of the
variables $\{\tX_{1j}\dvtx1\leq j \leq p\}$,
defined by
$\check{\varrho}_n(k) = \sup\{| P(A\cap B) - P(A)P(B)|/\sqrt{P(A)
P(B)} \dvtx A \in\F_1^m, B\in\F_{m+k}^p, 1\leq m\leq p-k\}
$. Here, $\F_a^b$ denotes the $\sigma$-field generated by $\{\tX_{1j} \dvtx a\leq j \leq b\}$, $1\leq a\leq b\leq p$, $\tX_{1j}
=X_{1\tau(j)}, $
and $\tau=\tau_n$ is the inverse of the permutation $\pe$
in (C.4)$_\al$.
\end{longlist}

Condition (C.4)$_\al$ says that up to a (possibly unknown) permutation
of the
coordinates, the components of the
$X_i$-vectors have a dependence structure
that is asymptotically similar to the one given by
$\rho_{\al}(\cdot)$. Note that
the sum $\sum_{k=0}^\infty|\rho_{\al}(k)|$ diverges
if and only if $\al\leq1$,
and therefore, we classify the dependence structure of the
$X_1j$'s as LRD or SRD according to $\al\leq1$ or $\al>1$,
respectively; cf. \citet{B94}. Condition (C.5)$_\al$ is
a decay condition on the
$\varrho$-mixing coefficient of the reordered
variables $\{\tX_{1j}\dvtx1\leq j\leq p\}$. Note that
by (C.4)$_\al$, the correlation coefficient
between $\tX_{1j}$ and $\tX_{1k}$ is $\rho_{\al}(j-k)(1+o(1))$,
and therefore, the re-ordered sequence $\{\tX_{1j}\dvtx1\leq j\leq p\}$
behaves approximately like
a stationary time series with the natural time-index $j$.
Thus, condition (C.5)$_\al$ specifies the
degree of dependence of the $X_{ij}$'s,
up to a permutation that need not be known to the
user.


\subsubsection{Results under short-range dependence}\label{sec3.3.1}
The following result shows that
the log-PEL ratio statistic
in the SRD case
converges to a Normal limit after a suitable centering and
after scaling by the ``standard factor'' $p^{1/2}$.

\begin{theorem}\label{th3.2}
Let conditions \textup{(C.1), (C.2)}
and \textup{(C.4)}$_\al$,
\textup{(C.5)}$_\al$ hold for some $\al> 1$, $s \ge6, r \ge12 $.
Let $\ka^2 = 2c_*^2\sum_{k=0}^\infty\rho_0(k)^2$.
Then, for $p = o(n^2)$,
%
\begin{equation}
p^{1/2} \bigl[ - \log R_n(\mu_0) - c_*
  \bigr] \raw^d N \bigl(0, \ka^2 \bigr).
\label{wk-dep}
\end{equation}
%
\end{theorem}

We now comment on Theorem~\ref{th3.2}. From the proof,
it follows that the distribution of the log-PEL
ratio statistic, for a given sample size, is close to
the sum of $p$ weakly dependent chi-squared random variables
with one degree of freedom. As a result, centering at 1
and scaling by $p^{1/2}$ yields a\vadjust{\goodbreak} nondegenerate Normal limit.
The effect of the weak dependence shows up in the
variance of the limiting Normal distribution, which
depends on the correlation structure of
the components of $X_{1}$. It is worth noting
that in the SRD case, one can use
Normal critical points with an estimated variance
to calibrate simultaneous tests of $p$ hypotheses
using the EL.

Theorem~\ref{th3.2} extends existing results on the EL
in more than one direction.
Hjort, McKeague and Van~Keilegom (\citeyear{H09}) and Chen et al. (\citeyear{C09}) proved a version of the result
(i.e., a Normal limit) for the standard log-EL
ratio statistic in increasing dimensions
with centering at 1 and scaling by $p^{1/2}$.
In comparison, Theorem~\ref{th3.2} relaxes the restriction on
the dimension $p$ of the parameter $\mu$, by allowing it
to grow faster than the sample size. This
should be compared with the best available
rate of $p=o(n^{1/2})$,
obtained by Chen et al. (\citeyear{C09}).
Further, Theorem~\ref{th3.2} covers a wide range of dependence structures
of the components of $X_{1j}$'s which are not covered by the earlier
results (e.g., here the minimum eigen-value of the covariance matrix of
$X_1$ need not be bounded away from zero). However, the most important
implication of Theorem~\ref{th3.2} is that under SRD,
the penalization step circumvents the limitation
of the standard EL which is known to break down beyond the threshold
$p\leq n/2$, as shown by \citet{Tsao04}.

\subsubsection{Results under long-range dependence}\label{sec3.3.2}
For $\al\in(0,1]$, the sum\break
$\sum_{k=1}^{\infty} \rho_0 (k)$ fails to converge absolutely,
and we refer to this as
the LRD case.
Sums of LRD random variables are known to have
either a Normal or a
non-Normal limit, depending on the value of $\al$.
The next result deals with the case where $\al$ can be
very small, and the limit law is non-Normal. Further, the scaling
also depends on the correlation decay parameter $\al$,
as shown by Theorem~\ref{th3.3}. Let $\Ga(\al) = \int_0^\infty t^{\al
-1}e^{-t}\,dt$
and ${\iota}=\sqrt{-1}$.

\begin{theorem}\label{th3.3}
Let conditions \textup{(C.1), (C.2)} and \textup{(C.4)}$_\al$,
\textup{(C.5)}$_\al$ hold for some $\al\in(0,\frac{1}{2})$, $s \ge6$
and $r \ge12$. If $p = o(n^{1/\al})$, then
%
\begin{equation}
p_n^{\al} \bigl[ - \log R_n(\mu_0)
- c_* \bigr] \raw^d W, \label{st-dep}
\end{equation}
where $W$ is defined in terms of a bivariate Wiener--It\^o
integral with respect to the random spectral measure $\zeta$
of the Gaussian white noise process
as
\[
W = c_* \bigl[2\Ga(\al) \bigr]^{-1}\int\frac{\exp({\iota}[x_1+x_2]) -
1}{{\iota}[x_1+x_2]}
|x_1x_2|^{(\al-1)/2} \,d\zeta(x_1) \,d
\zeta(x_2). 
\]
\end{theorem}

Theorem~\ref{th3.3} shows that
under very strong dependence (i.e.,
for small values of~$\al$) in the E-case,
the log-PEL ratio statistic, with the same centering but a different
scaling factor, has a nondegenerate limit distribution and the limit law
is \textit{non-Normal}. Further, the range $p$ for which the
result holds is $p= o(n^{1/\al})$, which is a decreasing
function of $\al$. Thus, the stronger the dependence
among the co-ordinates of $X_{1j}$'s, the larger is the
allowable growth rate of $p$ as a function of $n$
for the validity of the limit distribution.\vadjust{\goodbreak}

%

Theorems~\ref{th3.2} and~\ref{th3.3} exhaust the types of
limit laws for the log-PEL ratio statistic
in the E-case. However, in terms of
the rate of decay of the correlation function, these
leave out the case where $\al\in[1/2,1]$.
Although $\al\in[1/2,1]$ corresponds to LRD
in the traditional sense,
the centered and scaled
versions of the log-PEL ratio statistic continue to have
a Normal limit as shown by the
following result.
Curiously, the
scaling sequence as well as the growth rate of $p$ depend on
whether $\al=1/2$ or $\al\in(1/2,1]$.

\begin{theorem}\label{th3.4}
Let conditions \textup{(C.1), (C.2)} and \textup{(C.4)}$_\al$,
\textup{(C.5)}$_\al$ hold for some $1/2\leq\al\leq1$, $s \ge6, r \ge12 $.
\begin{longlist}[(ii)]
\item[(i)]
If $1/2< \al\leq1$ and $p = o(n^2)$, then
\eqref{wk-dep} holds.
\item[(ii)]
If $\al=1/2$ and $p=o([n/\log n]^2)$, then
%
\begin{equation}
[p\log p]^{1/2} \bigl[ - \log R_n(\mu_0) - c_*
\bigr] \raw^d N \bigl(0,c_*^2 \bigr). \label{bdr-dep}
\end{equation}
\end{longlist}
\end{theorem}

Thus, it follows from Theorem~\ref{th3.4} that the log-PEL ratio statistic is
asymptotically Normal for all $\al\geq1/2$, although the components
of $X_{i}$'s have LRD when $\al\in[1/2,1]$. The peculiar
behavior of the scaling sequence
at the boundary value $\al=1/2$ is essentially
determined by the growth rate of the series $\sum_{j=1}^p \rho_\al^2(j)$ as
$p\raw\infty$, which is asymptotically equivalent to $\log p$
for $\al=1/2$ but it is bounded for $\al>1/2$.

\begin{remark}\label{rem3.1} Proofs of Theorems~\ref{th3.2}--\ref{th3.4} show that
for any $p\raw\infty$,
$-\log R_n(\mu_0) \raw_p c_* \mbox{ as } \nti$,
that is, the log-PEL ratio statistic has a degenerate limit
under \textit{all} sub-cases of the E-case, for \textit{arbitrarily
large} $p$ as a function of $n$. However, for nondegenerate
limits, refined approximations to the difference
$[-\log R_n(\mu_0) - c_*]$ are needed. Here, we are able to show that
an approximation of the form
%
\begin{equation}
-\log R_n(\mu_0) - c_* = T_n +
E_n \label{E-app}
\end{equation}
holds for \textit{all} sub-cases of the E-case, where $T_n$
is a centered sum and where $E_n$ is an error
term, roughly of the order of $O_p(n^{-1})$. Further,
$T_n$ has a nondegenerate limit up to a suitable scaling,
as a function of $p$, depending on the dependence
structure of the $X_{1j}$'s. The bounds on the growth rate
of $p$ in the different sub-cases of the E-case are then
determined by the requirement that the scaled error\vadjust{\goodbreak}
term be asymptotically negligible. For example, in the
SRD case, $p^{1/2} T_n\raw^d N(0, \ka^2)$ and hence,
$p^{1/2}E_n\raw_p 0$ if and only if $p^{1/2}/n \raw0$,
which is equivalent to the bound $p = o(n^{2})$. Similar
considerations lead to the respective upper bounds in
the other sub-cases of the E-case.
\end{remark}

\begin{remark}\label{rem3.2} It is worth pointing out that
the PEL can be used for
constructing {``conservative''} large sample
simultaneous tests of the $p$ hypotheses
$H_0\dvtx\mu= \mu_0$ for\vadjust{\goodbreak} \textit{arbitrarily}
large $p$ in the E-case. Indeed, for $p$ growing faster than
the upper bounds given in Theorems~\ref{th3.2}--\ref{th3.4}, a
conservative large sample simultaneous test of
$H_0\dvtx\mu= \mu_0$ rejects $H_0$ if
$|c_*+\log R(\mu_0)| > n^{-1}\log n$. Note that by
\eqref{E-app}, this test attains the \textit{ideal}
level $0$ asymptotically.
\end{remark}

\section{A subsampling based calibration}\label{sec4}
In this section, we describe a nonparametric
calibration method based on subsampling to
approximate the quantiles of the nondegenerate
limit laws in both E- and NE-cases, which typically
involve unknown population parameters
and hence, cannot be used directly in practice.
Let $\cX_n(I) = \{X_i \dvtx i\in I\}$ be a
subset of $\{X_1,\ldots,X_n\}$ where $I\subset\{1,\ldots,
n\}$ is of size $m$ and where $1< m < n$ (specific conditions
on $m$ are given below).
On each $\cX_n(I)$,
we employ the PEL method and obtain a version of the PEL ratio
statistic $R^*_m(\mu;I)$, by replacing $n$ with $m$ and
$X_1,\ldots, X_n$ by $\cX_n(I)$ in the definitions
$R_n(\mu)$. First
consider the NE case. Here,
the subsampling estimator of the distribution function $G_{n}^{\ne}
(\cdot) \equiv P(-\log R_n(\mu_0) \leq\cdot)$
under the null hypothesis $H_0\dvtx\mu= \mu_0$
is given by
\[
\hG_n^{\ne}(x) = |\cI_n|^{-1} \sum
_{I\in\cI_n} \ind \bigl( -\log R^*_m(
\mu_0;I) \leq x \bigr),  \qquad  x\in\bbr,
\]
where $\cI_n$ is a collection of subsets of
$\{1,\ldots,n\}$ of size $m$ and where $|A|$ denotes
the size of a set $A$. All possible subsets
of size $m$ cannot be used mainly due to
the sheer number of such sets, and hence,
only a small fraction of these subsets
are used to compute $\hG_n^{\ne}(\cdot)$
in practice. In view of the block resampling
methods for time series data,
here we shall take $\cI_n$ to be the collection of all overlapping
blocks (subsets) of size $m$ contained in $\{1,\ldots,n\}$. Then,
we have the following result 
in
the NE case.

\begin{theorem}\label{th4.1} Suppose that the conditions of Theorem~\ref{th3.1}
hold and that
%
\begin{equation}
m^{-1} + m/n = o(1) \qquad\mbox{as } \nti. \label{m-cond}
\end{equation}
Then,
$
\sup_{x\in\bbr}  | \hG_n^{\ne}(x) -G_n^{\ne}(x) | \raw_p 0
\mbox{ as } \nti.$
\end{theorem}

Next consider the E-case. Note that for $\al= 1/2$, the limit of the
log-PEL ratio statistic is $N(0,1)$, which is distribution
free.
One can carry out a simple test [cf. \citet{B94}] to ascertain\vadjust{\goodbreak}
if ``$H\dvtx\al= 1/2$'' is true and then use the limit
distribution directly to conduct the PEL test of the
simultaneous $p$ hypotheses $H_0\dvtx\mu=\mu_0$
using the $N(0,1)$ critical points, without the
need for an alternative calibration. As a result, we concentrate
on the values of $\al\neq1/2$ in the E-case.
Let $\hal_n$ be an estimator of the
correlation parameter $\al$; cf. Remark~\ref{rem4.1} below.
Let $R_m^*(\mu_0;I)$ denote the PEL ratio statistic
based on the subsample $\cX(I)$ under $\mu_0$, and define
\[
V^*_{m}(I) = \hb_n \bigl[ -\log R_m^*(
\mu_0;I) - c_* 
\bigr],  \qquad  I\in\cI_n,
\]
where $ \hb_n = p^{\hal_n\wedge{1}/{2}}$.
Then, a subsampling estimator of the distribution
of $V_n\equiv\hb_n[-\log R_n(\mu_0) - c_*] $ is given by
$\hG_{n,\al}(x)= |\cI_n|^{-1} \sum_{I\in\cI_n}
\ind( -\log V^*_{m}(I) \leq x),   x\in\bbr
$ and
we have the following results.


\begin{theorem}\label{th4.2}Suppose that there exists a $c_0\in
\bbr$
such that
%
\begin{equation}
(\log p) [\hal-\al] \raw_p c_0 \qquad\mbox{as }\nti.
\label{al-cond}
\end{equation}\vspace*{-12pt}
\begin{longlist}[(ii)]
\item[(i)]
For $\al\in(1,\infty)$, let the conditions of Theorem~\ref{th3.2}
and for $\al\in(1/2,1]$, let those of Theorem~\ref{th3.4}\textup{(i)}
hold. If $p/m^{2} + m/n = o(1)$, then
%
\begin{equation}
\sup_{x\in\bbr} \bigl| \hG_{n,\al}(x) - P (V_n \leq x ) \bigr|
\raw_p 0 \qquad{\mbox{as } } \nti. \label{ss-w-d}
\end{equation}
%
%
\item[(ii)]
If the conditions of Theorem~\ref{th3.3}
hold for some $\al\in(0,1/2)$
and $p^\al/m + m/n = o(1)$, then \eqref{ss-w-d} holds.
\end{longlist}
\end{theorem}

Theorem~\ref{th4.2} shows that for {both} Normal and non-Normal
limit laws under the E-case, the
subsampling method provides a valid approximation to
the distribution of the log-PEL ratio statistic.
Hence, one can use the quantiles of the subsampling
estimators to calibrate simultaneous
tests on $\mu$ in a unified manner. This is specially important in
the case of non-Gaussian limit laws for which the quantiles
are difficult to derive. However, for $\al> 1/2$,
the limit distribution is Gaussian, and an
alternative approximation can be generated by using a Normal
distribution with an \textit{estimated} variance.
Indeed,
the latter
may be preferable to subsampling
from the computational point of view.

\begin{remark}\label{rem4.1}
In practice, the value of $\al$ is not known and
must be estimated. First consider the
case where the permutation $\iota_n(\cdot)$ in
(C.4)$_\al$ is known. Then, we are essentially dealing with
$n$ i.i.d. copies of a time series of length $p$
as observations. By using the $n$ replicates of
the time series, it is easy to modify
standard estimators of $\al$ based on
a single time series [cf. \citet{B94}] to construct
an estimator $\hal$ of $\al$ satisfying
$\hal- \al= o_p(n^{-1/2}[\log p]^{-1}) \mbox{ as } \nti
$,
which clearly satisfies \eqref{al-cond} with $c_0=0$.

Next consider the case where the permutations
$\iota_n(\cdot)$ are unknown. In this case, it is \textit{not}
possible to identify pairs $(j,l)$, $1\leq j,l\leq p$ that correspond to
the lag-$k$ correlation $\rho_\al(k)$. However, it is still possible to
construct estimators of $\al$ that satisfy \eqref{al-cond}. Define
%
\begin{equation}
\hal= - (\log p)^{-1} \log \Biggl( e_n + n^{-1}
\sum_{i=1}^n \Biggl\{p^{-1}\sum
_{j=1}^p [X_{ij} -
\bX_{nj}]\delta_j^{1/2} \Biggr\}^2
\Biggr), \label{p-inv-hal}
\end{equation}
where $e_n = \prod_{j=1}^p \ind(s_{nj}=0)$ and $\delta_j$
(and $\bX_{nj}$ and $s_{nj}$)
are as in \eqref{de-j}. Note that $\hal$ is invariant under
permutations of the components of $X_i$'s and also under component-wise\vadjust{\goodbreak}
location and scale transformations. In Section~\ref{sec6}, we show that
$\hal$ satisfies~\eqref{al-cond} under the conditions of Theorem~\ref{th4.2},
even when $\iota_n(\cdot)$ is unknown. In the same spirit,
we may use the following
estimator of the limiting variance $\ka^2$ in
Theorem~\ref{th3.2} in the case where $\iota_n(\cdot)$ is unknown:
%
\begin{equation}
\hat{\ka}^2 = 2c_*^2 \sum_{j=2}^{p-1}
\hat{c}(1,j)^2 \ind \bigl(\bigl|\hat{c}(1,j)\bigr|>2n^{-1/2}\log n
\bigr), \label{p-inv-hka}
\end{equation}
where $\hat{c}(j,k) = n^{-1}\sum_{i=1}^n (X_{ij}-\bX_{nj})(X_{kn}-\bX_{nk})
\delta_j \delta_k$, $1\leq j,k\leq p$. Consistency of $\hat{\ka}^2$
holds under mild moment conditions; cf. Section~\ref{sec6}.
\end{remark}

\begin{remark}\label{rem4.2} An important factor that impacts the
accuracy of the subsampling
method is
the choice of the subsample size $m$.
Note that
%
\begin{equation}
m= C [np]^{\al_0/(1+\al_0)} \label{m-ch}
\end{equation}
satisfies the requirements of 
Theorem~\ref{th4.2},
where $\al_0=\min\{\al,1/2\}$. At this point,
we do not know the order of the optimal $m$ for the
different cases considered here.
In the next section, we address this through a numerical
study and explore the effects of different choices of $m$
on the performance of the PEL method.
\end{remark}


\section{Numerical study}\label{sec5}
We assess finite sample performance of the PEL
method by simulation in a variety of settings. We considered
different combinations of the sample size $n$ and the dimension $p$, with
$p \approx2n^{1/3}, n/2, 2n$ and $n= 40$ and $200$.
The testing problem we considered is $H_0\dvtx \mu= 0$,
although any other value of $\mu$ may be used in $H_0$,
as the PEL criterion is location invariant.
We generated i.i.d. random $p$-vectors $X_1,\ldots,X_n$
where the $p$ coordinates
of $X_i$'s had one of the three different types
of dependence structures,
namely: (i)
non-Ergodic, (ii) LRD
and (iii) SRD, as follows.

\subsection{Algorithms for generating the data}\label{sec5.1}
\subsubsection{The nonergodic case}\label{sec5.1.1}
\begin{longlist}[(1)]
\item[(1)] Consider the basis functions in $L^2[0,1]$ given by
$\phi_{j}(t)=\sin{(2\pi jt)}/\sqrt{2}$ {for}  $j=1,2,\ldots,15$
{and}  $\phi_{j}(t)=\cos{(2\pi(j-15)t)}/\sqrt{2}$,  {for}
$j=16,2,\ldots,30 $ and $\phi_{0}(t)\equiv1$.\vadjust{\goodbreak}
\item[(2)]
Generate $Z_j \sim\rm{N}(0,1)$ i.i.d. and let
$\lambda_j = ( \exp{\{1\}}+1  ) $ for all $j$.
\item[(3)] Define $X_{1j}= \sigma_j \cdot W(j/p),
 j=1,2,\ldots, p$, where
$W(\cdot)=\sum_{j=0}^{30}Z_j \phi_j({\cdot}) \lambda_j$
and where $\si_j\equiv\si_{jn}$ are scalars in $(0,\infty)$.
\end{longlist}
Then, $X_1= \{ X_{11},X_{12},\ldots , X_{1p}  \} $ is
a nonergodic series. Replicates of $X_1$ yield
$X_1,\ldots,X_n$ in the NE-case.
\subsubsection{Long-range dependence}\label{sec5.1.2}
For the LRD case, we follow a setup similar to that used in
Hall, Jing and Lahiri (\citeyear{HJL98}). We generate stationary increments
of a self-similar process with self-similarity parameter\vadjust{\goodbreak}
(or Hurst constant)
$H = \frac{1}{2}(2 - \alpha)\in(1/2,1)$ for $\al\in(0, 1)$.
The algorithm is as follows:
\begin{longlist}[(1)]
\item[(1)] Generate a random sample $\textbf{Z}_{p0} = \{
Z_{10},\ldots,
Z_{p0}\}$
from $N(0,1)$.
\item[(2)]
Define ${ \bf Z}_p \equiv U^{T}{ \bf Z}_{p0}$, where $U$ is
obtained by Cholesky factorization of $R$ into $R=U^{T}U$
and where $R=((r_{ij}))$ with $r_{ij}=\rho_\al(|i-j|)$,
and
%
\begin{equation}
\rho_\al(k) = \tfrac{1}{2} \bigl\{ (k+1)^{2H} +
(k-1)^{2H} - 2k^{2H} \bigr\}, \qquad  k\geq1,
\end{equation}
and $\rho_\al(0)=1$. Note that $\rho_\al(k)
\sim Ck^{-\alpha}  \rm{as}  k \rightarrow\infty$.

\end{longlist}
Replicates of $\textbf{Z}_p$ give the variables $X_1,\ldots,X_n$ in the
LRD case.

For the simulation study here, we considered
the NE-case ($\al=0$) where the data were generated
by the algorithm in Section~\ref{sec5.1.1} and the LRD cases
$\al= 0.1$ and $0.8$ based on
the algorithm of Section~\ref{sec5.1.2}. For the SRD case
($\al=\infty$), $X_1$ was generated by an
$\operatorname{ARMA}(2,3)$ process with $N(0,1)$ error variables
and parameter vector $(-0.4, 0.1; 0.3,0.5,0.1)$.

\subsection{Choice of the subsample size}\label{sec5.2}
We also considered different choices of the subsample size $m$
in order to get some insight into its effects on the accuracy of the
subsampling calibration. Note that the feasible choices of
the subsample size depend on the relative growth rates of
both $n$ and $p$ as well as on the
strength of dependence, here quantified by $\al$. For each pair
$(p,n)$, we considered three choices of the subsample size $m$
(denoted by the generic symbols $m_1,m_2,m_3$), depending on the dependence
structure. Specifically, for the SRD case ($\al=\infty$), we set
\[
m_i = c_i^0 \cdot[np]^{1/3},\qquad
  i=1,2,3,
\]
where $c_1^0=0.5$, $c_2^0 =1$ and $c_3^0 = 2$. Note that in this case,
the random variables in $X_1,\ldots,X_n$ form a series of length $np$ and
are weakly dependent. Further, the target parameter for the subsampling
method in the SRD case (and also in the LRD case with
$\al>1/2$) is the variance of the limiting
Normal distribution. Hence, in view of the well-known results
on optimal block length (for variance estimation)
[cf. \citet{HHJ99}, \citet{L03}],
the above choices of the $m_i$'s
are reasonable.

\begin{table}
\caption{Empirical levels of significance $\hat{a}$
for the subsampling based PEL, with sample size $n=200$
at $0.05$ significance level and $c_*=1$. Here we have
reported $|0.05 - \hat{a}|$}\label{tabc=1}
\begin{tabular*}{\textwidth}{@{\extracolsep{\fill}}lccccccccc@{}}
\hline
& \multicolumn{3}{c}{$\bolds{p_1=20}$} & \multicolumn{3}{c}{$\bolds{p_2=100}$}
& \multicolumn{3}{c@{}}{$\bolds{p_3=400}$} \\[-6pt]
& \multicolumn{3}{c}{\hrulefill} & \multicolumn{3}{c}{\hrulefill}
& \multicolumn{3}{c@{}}{\hrulefill} \\
$\bolds{\alpha}$ & $\bolds{m_1}$ & $\bolds{m_2}$ &
$\bolds{m_3}$ & $\bolds{m_1}$ & $\bolds{m_2}$ &
$\bolds{m_3}$ & $\bolds{m_1}$ & $\bolds{m_2}$ &
\multicolumn{1}{c@{}}{$\bolds{m_3}$}\\
\hline
$0.0 $& 0.0092 &0.0063 &0.0190 &0.0090 &0.0114 &0.0162 &0.0089 &0.0103
&0.0091 \\
$0.1$& 0.0099 &0.0081 &0.0171 &0.0045 &0.0069 &0.0129 &0.0149 &0.0134 &0.0094 \\
$0.8$& 0.0079 &0.0031 &0.0061 &0.0086 &0.0042 &0.0081 &0.0101 &0.0099 &0.0190\\
$\infty$& 0.0059 &0.0091 &0.0010 &0.0020 &0.0104 &0.0039 &0.0091 &0.0159
&0.0078 \\
\hline
\end{tabular*}
\end{table}

\begin{table}[b]
\caption{Empirical levels of significance $\hat{a}$
for the subsampling based PEL, with sample size $n=200$
at $0.1$ significance level and $c_*=1$. Here we have
reported $|0.1 - \hat{a}|$}\label{tabn200}
\begin{tabular*}{\textwidth}{@{\extracolsep{\fill}}lccccccccc@{}}
\hline
& \multicolumn{3}{c} {$\bolds{p_1=20}$} & \multicolumn{3}{c} {$\bolds{p_2=100}$}
& \multicolumn{3}{c@{}} {$\bolds{p_3=400}$} \\[-6pt]
& \multicolumn{3}{c} {\hrulefill} & \multicolumn{3}{c}{\hrulefill}
& \multicolumn{3}{c@{}}{\hrulefill} \\
$\bolds{\alpha}$ & $\bolds{m_1}$ & $\bolds{m_2}$ &
$\bolds{m_3}$ & $\bolds{m_1}$ & $\bolds{m_2}$ &
$\bolds{m_3}$ & $\bolds{m_1}$ & $\bolds{m_2}$ &
$\bolds{m_3}$\\
\hline
$0.0 $&0.170 &0.020 &0.075 &0.221 &0.142 &0.090 &0.075 &0.152 &0.227\\
$0.1$& 0.030& 0.033 & 0.033 & 0.012 & 0.005 &0.011 &0.112 &0.133 &0.066\\
$0.8$&0.011 &0.045 &0.087 & 0.123& 0.082&0.018 &0.108 &0.027 & 0.069\\
$\infty$&0.138 & 0.135& 0.065&0.050 &0.011 &0.003 &0.048 &0.054 &0.026\\
\hline
\end{tabular*}
\end{table}
%

Next consider the case $\al=0.1$ under LRD, where the limit
distribution is
non-Normal. From the proofs of Theorems~\ref{th3.2} and~\ref{th4.2},
it follows that the prescription for $m$ in \eqref{m-ch}
attempts to balance the bias of the subsampling
approximation to the limit distribution and its variance.
However, for a very small value of $\al$, a direct application
of \eqref{m-ch} leads
to a very small fractional exponent of $np$, which may be too
small in practice.
In such situations, particularly where $p$ is not very
large and the LRD exponent $\al$ is small, we use the threshold
$n^{1/3}$ and set
\[
m_i\equiv m_i (\al) = c_i^0
\cdot\max \bigl\{[np]^{\al/(1+\al)}, n^{1/3} \bigr\}, \qquad  i=1,2,3,
\]
where $c_i^0$'s are as before. The rationale behind this modification
is that for $p$ small, we simply treat $X_1,\ldots, X_n$ as a
weakly dependent
multivariate time series
and again employ the known results on the
optimal block size.\vadjust{\goodbreak} 

%
%
%

Finally, consider the NE-case, $\al=0$. Note that for $\al=0$,
$p$ can grow at an arbitrary rate with the sample size $n$
for the validity of Theorems~\ref{th3.1} and~\ref{th4.1}. Hence,
in this case, our choice of $m$ depends only
on the sample size. We consider the ``canonical''
choice $m_1=n^{1/3}$ as well as the larger values
$m_2= n^{1/2}$ and $m_3= 2n^{1/2}$ to explore the
effects of a larger subsample size on the accuracy
of the subsampling calibration method.

\subsection{Results}\label{sec5.3}
\subsubsection{Levels of significance in simultaneous tests}\label{sec5.3.1}
Here we consider finite sample accuracy of the proposed
PEL method for simultaneous testing of
the $p$ hypotheses
%
\begin{equation}
H_0\dvtx\mu=0 \quad\mbox{vs.}\quad  \mu\neq0 \label{hyp-1}
\end{equation}
at the levels of significance
$a=0.1, 0.05$.
The correlation parameter $\alpha$ for the
subsampling based calibration was estimated by
averaging the
Taqqu, Teverovsky and Willinger (\citeyear{Taq95}) estimator of
the Hurst parameter ($H$) from each of the
$p$-time series and by setting
$\hat\alpha= (2-2 \hat H) $. Further,
we have used the interior-point method
[cf. \citet{Wr97}]
as a fast optimization tool for computing the PEL ratio statistic,
which can handle high-dimensional optimization
problems efficiently.
Tables~\ref{tabc=1} and~\ref{tabn200} report the\vadjust{\goodbreak}
attained levels of significance based on $500$
simulation runs and $n=200$ for the target significance levels
of $0.05$ and $0.10$, respectively, for
different values of $p$, $m$, and $\al$.

\begin{table}[b]
\caption{Power of the proposed PEL, with sample size $n=200$ at $0.1$
significance level and $c_*=1$}\label{tabpower}
\begin{tabular*}{\textwidth}{@{\extracolsep{\fill}}lccccccccc@{}}
\hline
& \multicolumn{3}{c}{$\bolds{p_1=20}$} & \multicolumn{3}{c}{$\bolds{p_2=100}$}
& \multicolumn{3}{c@{}}{$\bolds{p_3=400}$} \\[-6pt]
& \multicolumn{3}{c}{\hrulefill} & \multicolumn{3}{c}{\hrulefill}
& \multicolumn{3}{c@{}}{\hrulefill} \\
$\bolds{\alpha}$ & $\bolds{m_1}$ & $\bolds{m_2}$ &
$\bolds{m_3}$ & $\bolds{m_1}$ & $\bolds{m_2}$ &
$\bolds{m_3}$ & $\bolds{m_1}$ & $\bolds{m_2}$ &
\multicolumn{1}{c@{}}{$\bolds{m_3}$}\\
\hline
$0.0 $& 0.569 &0.681 &0.929 &0.515 &0.643 &0.791 &0.87\phantom{0} &0.834 &0.766 \\
$0.1$& 0.66\phantom{0} &0.71\phantom{0} &0.85\phantom{0} &0.569 &0.903 &0.676 &0.794 &0.8\phantom{00} &0.868 \\
$0.8$& 0.515 &0.883 &0.688 &0.75\phantom{0} &0.997 &0.488 &0.70\phantom{0} &0.87\phantom{0} &0.90\phantom{0}\\
$\infty$& 0.510 & 0.622 & 0.870 & 0.739 & 0.778 & 0.996 & 0.802 &
0.790 & 0.939
\\
\hline
\end{tabular*}
\end{table}

From the tables, it follows that the PEL does a reasonable job of
simultaneous testing of $p$ hypotheses
for all $4$ cases of dependence, for appropriately
chosen subsample size.
Comparing the attained level of significance, it is clear that
the best choice of the subsample size
critically depends on the relative
sizes of $n$ and~$p$, and more importantly, on the
type of dependence among the components of $X_1$.
Further, rather surprisingly,
the PEL tests at the level of significance $0.05$
turned out to be more accurate (on an absolute scale)
than at the level $0.1$, for the subsample sizes considered here.

%


We also considered the effect of
the penalty parameter $\lambda_n=c_*n/p$
on the performance of the PEL test. In the supplementary material \citet{LM12} (hereafter referred to as [LM]), we report
the empirical levels of significance of the PEL test
for $n=200$ and the target level $0.1$ for two other
choices of the constant $c_*$, namely, $c_*=0.5$
and $c_*= 2.0$. The results for the choice $c_*=2$
are qualitatively similar to those reported in Table
\ref{tabn200} (with $c_*=1$); in comparison, the accuracy for the
case $c_*=0.5$ appears to be slightly better than
the $c_*=1$ case. A similar pattern was
observed for the $0.05$ level of significance.
We also considered the accuracy of the empirical
significance levels of the PEL tests at a
relatively smaller sample size $n=40$, for
$c_*=1$ and $\al=0.1$; cf. [LM]. The PEL has a
reasonable performance
even at this low sample size;
see [LM] for details.

\subsubsection{Finite sample power properties}\label{sec5.3.2}
To get some idea about the power properties of the
PEL tests, we computed the probability
of Type II error for a level $0.1$ PEL test
with $n=200$ and $c_*=1$
under the alternative
%
$\mu=\mu_1$ where the first $p/2$ components of $\mu_1$ were
equal to $1$ and the rest were $0$.
Table~\ref{tabpower} gives the power of the PEL test
at level $0.1$ under $\mu=\mu_1$
for
different combinations of $p$, $\al$ and $m$.
From Table~\ref{tabpower}, it appears that the power
can be reasonably high for a suitable choice
of the subsample size, although the maximum value critically
depends on the dimension $p$ of the parameters
and the strength\vadjust{\goodbreak} of dependence $\al$. In particular, the PEL
attains a higher (maximum) power under
weaker dependence ($\al=0.8,\infty$) than
under strong dependence ($\al=0,0.1$).

\subsubsection{Comparison with Normal calibration}\label{sec5.3.3}
Note that for $\al>1/2$, the limit distribution of the
logarithm of the PEL ratio statistic is Normal and
therefore, one can use the limiting Normal distribution
with an estimated variance to conduct the PEL
test. In this section, we compare the performance of the
subsampling-based calibration with the Normal
distribution-based calibration.
To estimate the asymptotic variance
$
\kappa^2=2c_*^2\sum_{k=0}^\infty\rho_0(k)^2
$,
we first estimate $\rho_0(k)$ using the
sample auto-covariance at lag-$k$ based on the components of
individual $X_i$'s and then average them to get an estimate
${\hat{\rho}}_n(k)$ of $\rho_0(k)$ for $k=1,\ldots,K$
where $K=\min\{p/2, p^{1/2}\}$.
Since $\ka^2$
involves the squares of $\rho_0(k)$, the plug-in estimator
is positive (with probability $1$).
Tables~\ref{tabES40} and~\ref{tabES200} compare the best
performance of the subsampling based PEL with the
Normal, calibration-based
PEL for $n=40$ and $n=200$, respectively.

\begin{table}
\caption{Comparison of the subsampling (SS) and
Normal (G) calibrations for $n=40$}\label{tabES40}
\begin{tabular*}{\textwidth}{@{\extracolsep{\fill}}lcccccc@{}}
\hline
& \multicolumn{2}{c} {$\bolds{p_1=7}$} & \multicolumn{2}{c} {$\bolds{p_2=20}$}
& \multicolumn{2}{c@{}} {$\bolds{p_3=80}$} \\[-6pt]
& \multicolumn{2}{c} {\hrulefill} & \multicolumn{2}{c} {\hrulefill}
& \multicolumn{2}{c@{}} {\hrulefill} \\
$\bolds{\alpha}$ & \textbf{G} & \textbf{SS}& \textbf{G}
& \textbf{SS} & \textbf{G} & \textbf{SS}\\
\hline
$0.8 $ & \textbf{0.122} & 0.151 & 0.132 & \textbf{0.080} & 0.136 &
\textbf{0.081}\\
$\infty$&\textbf{0.098} & 0.092 &0.030 &\textbf{0.111} &0.076 &\textbf{0.093} \\
\hline
\end{tabular*}
\end{table}

Tables~\ref{tabES40} and~\ref{tabES200} show that,
except for the small values of $p$, the subsampling-based PEL method
has a
better accuracy
(marked as bold) than the Normal, calibration-based PEL. However, the
computational burden associated
with the subsampling method is typically larger than the Normal-based PEL.
%

%
\begin{table}[b]
\caption{Comparison of the subsampling (SS) and Normal (G) calibrations
for $n=200$}\label{tabES200}
\begin{tabular*}{\textwidth}{@{\extracolsep{\fill}}lcccccc@{}}
\hline
& \multicolumn{2}{c} {$\bolds{p_1=7}$} & \multicolumn{2}{c} {$\bolds{p_2=20}$}
& \multicolumn{2}{c@{}} {$\bolds{p_3=80}$} \\[-6pt]
& \multicolumn{2}{c} {\hrulefill} & \multicolumn{2}{c} {\hrulefill}
& \multicolumn{2}{c@{}} {\hrulefill} \\
$\bolds{\alpha}$ & \textbf{G} & \textbf{SS}& \textbf{G}
& \textbf{SS} & \textbf{G} & \textbf{SS}\\
\hline
$0.8 $ & 0.150 & \textbf{0.113} & 0.141 & \textbf{0.082} & 0.174 & \textbf{0.127}\\
$\infty$ &\textbf{0.111} & 0.165 & 0.048 & \textbf{0.103} & 0.238 & \textbf{0.126}\\
%
\hline
\end{tabular*}
\end{table}




\section{Proofs}\label{sec6}
Note that for each $n\geq1$, the PEL likelihood function
in \eqref{pel} is invariant with respect to (i) component-wise scaling
and (ii) permutation of the $p$ components. Hence, all through this
section, without loss of generality (w.l.g.), we set the
component-wise variance $\si_{nj}^2=1$ and set the permutation
$\pe(j) = j$ for all $1\leq j\leq p$ and $n\geq1$.
Let $C, C(\cdot)$ denote\vadjust{\goodbreak} generic constants that
depend only on their arguments (if any), but not on $n$.
Unless otherwise specified,
dependence on
(limiting) population quantities [such as $\rho_0(\cdot)$,
mixing coefficients, etc.] are dropped to simplify notation,
and limits in all order symbols are taken
by letting $\nti$. For $x\in\bbr$, let $\lfloor x\rfloor$
denote the largest integer not exceeding $x$
and let $x_+ = \max\{x,0\}$.

\subsection{Limit distribution in the nonergodic case}\label{sec6.1}
%
\begin{lemma}
\label{lem-1}
For each $n \ge1$, let \{$ \Yij= \Y_{ijn}, 1 \le j \le p $\},
$i=1,2,\ldots, n$ be a collection of $p=p_n$-dimensional random vectors
with $\E Y_{1j} =0 $ and $\E Y^2_{1j} =\si^2_{nj} \in(0,\infty)$. Let
$\delta_j \equiv\delta_{nj} = s^{-2}_{nj} \ind(s_{nj} \neq0) $ where
$s^2_{nj}={n^{-1}\sum_{i=1}^n}(Y_{ij}-\bY_{jn})^2$ and $\bY_{nj}={n^{-1}\sum_{i=1}^n}Y_{ij}$.
Also, let
$\Zij= |\Yij|-\E|\Yij|,$ $\bZ_{jn}={n^{-1}\sum_{i=1}^n}\Zij$,
$W_{i}(j,l) =
Y_{ij}Y_{il} - E Y_{ij}Y_{il}$, $\bW_{n}(j,l) = n^{-1}\sum_{i=1}^nW_{i}(j,l)$
and $D_{jn} =\{y\in\rl\dvtx \break P(Y_{1j}=y)>0\}$.
Suppose that
\textup{(L.1)}
$\max \lbrace\E ( {\si^2_{nj}}{\delta_{j}}  )^4
\dvtx1\le j \le p  \rbrace= O(1)$; \textup{(L.2)}
$ \max \lbrace\E|\Yij|^8 \dvtx1\le j \le p  \rbrace<
\C$,
and \textup{(L.3)} $\limsup_{\nti} $
$\max\{P(Y_{1j}=y)\dvtx y\in D_{jn}, j=1,\ldots, p\} <1$.
Then:
\begin{longlist}[(a)]
\item[(a)]
For $k=1,2,3$,
$ \sum_{j=1}^{p}\delta_j\bY^{2k}_{jn} =O_p(n^{-k}p)$.
\item[(b)]
$ \sum_{j=1}^{p}\delta_j[\bZ_{jn}^2+ \bW^2_{jn}] = O_p(n^{-1}p)$.
\item[(c)]
$ \sum_{j=1}^{p}\sum_{l=1}^{p}(\delta_j +1) \delta_l \bW^2_n(j,l) =
O_p (n^{-1}p^2 )$.
\item[(d)]
For $r=1,2$,
$\max_{1 \le j \le n}  \sum_{j=1}^{p}\delta_j|\Y_{1j}|^r = O( [pn]^{1/4})$.
\item[(e)]
$ \sum_{j=1}^{p}(\delta_j - 1)^2 = O_p (n^{-1}p )$.
\end{longlist}
\end{lemma}

\begin{pf}By
replacing $\Yij$'s with $\Yij/\si_{nj}$ for all $i,j$, w.l.g., we
can assume that $\si_{nj}=1$ for all $j,n$.
First consider part (a), $k=3$; the proofs of $k=1,2$ are similar.
By repeated use of H\"{o}lder's inequality,
\begin{eqnarray*}
\E\sum_{j=1}^p \delta_j |
\bY_{nj}|^6 \le \E \Biggl( \sum
_{j=1}^p \delta_j^4
\Biggr)^{1/4} \Biggl( \sum_{j=1}^p
\bY_{nj}^8 \Biggr)^{3/4} \le \Biggl(\E\sum
_{j=1}^p \delta_j^4
\Biggr)^{1/4} \Biggl(\E\sum_{j=1}^p
\bY_{nj}^8 \Biggr)^{3/4},
\end{eqnarray*}
which is $O(pn^{-3})$, by (L.1), (L.2). This proves (a).
Parts (b) and (c) follow by similar arguments.
As for part (d), note that
(for $r=2$)
\begin{eqnarray*}
\E \Biggl( \dm\max_{1 \le i \le n} \sum_{j=1}^p
\delta_j Y_{ij}^2 \Biggr) 
&\le& \Biggl[ n\E \Biggl( \sum
_{j=1}^p \delta_j Y_{1j}^2
\Biggr)^4 \Biggr]^{1/4}
\\
&\le& n^{1/4} \Biggl[ \E \Biggl(
\sum_{j=1}^p \delta_j^4
\Biggr)^{1/4} \Biggl( \E\sum_{j=1}^p
|Y_{1j}|^{8/3} \Biggr)^{3/4} \Biggr]^{1/4}
\le C n^{1/4} p^{1/4}.
\end{eqnarray*}
%
%
Finally consider part (e). Note that for $j=1,2,\ldots, p$
\begin{eqnarray*}
\delta_j - 1 &=&\frac{1- s_{nj}^{2} }{ s_{nj}^{2}} \ind(s_{nj} \neq0) -
\ind(s_{nj} = 0)
\\
&=& \delta_j \Biggl[ \bY_{nj}^2  -
{n^{-1}\sum_{i=1}^n} \bigl(
\Yij^2 - 1 \bigr) \Biggr] - \ind(s_{nj} = 0)
\\
&= & \delta_j \bY_{nj}^2 -
\delta_j \olW_n(j,j) - \ind(s_{nj} = 0) ,
\end{eqnarray*}
so that
%
\begin{equation}
\label{delta-ex} \delta_j = \ind(s_{nj}\neq0) +
\delta_j \bY_{nj}^2 - \delta_j
\olW_n(j,j).
\end{equation}
Part (e) can now be proved using (L.1)--(L.3)
and the Cauchy--Schwarz
inequality. We omit the details to save space.
\end{pf}


\begin{lemma}
\label{lem-2} Under the conditions of Theorem~\ref{th3.1},
\[
- \log{R}_n(\mu_0) = n\gamma_n \biggl(
\frac{\bY_{n1}}{\si_{n1}}, \ldots,\frac{\bY_{np}}{\si_{np}} \biggr) ( \I_p + 2
\gamma_n \A_n )^{-1} \biggl( \frac{\bY_{n1}}{\si_{n1}} ,
\ldots, \frac{\bY_{np}}{\si_{np}} \biggr)' + o_p(1),
\]
where $\bY_{nj}={n^{-1}\sum_{i=1}^n}Y_{ij}$, $\Yij= \Xij- \mu_j$,
$1 \le j \le p, 1 \le i \le n$ and
$\A_n = ((\rho_n(i-j)))_{p\times p}$.
\end{lemma}

\begin{pf}
W.l.g., let $\si_{nj}=1$ for all $1 \le j \le p, 1
\le i \le n$. Note that
%
\begin{equation}
- \log{R}_n(\mu_0) = \min \bigl\lbrace f(
\pi_1,\ldots, \pi_n) \dvtx(\pi_1,\ldots,
\pi_n)' \in\Pi_n \bigr\rbrace,
\end{equation}
where $f(\pi_1,\ldots, \pi_n) = - \sum_{i=1}^n\log(n \pi_i)
+ \lambda_n  \sum_{j=1}^{p}\delta_j  ( \sum_{i=1}^n\pi_i \Yij
)^2$.
Since $f(\cdot)$ is strictly convex in $\pi_1,\ldots, \pi_n$ over a
closed convex set $\Pi_n \subset\bbr^n$, it has a unique minimizer
in $\Pi_n$. [The maximum of $f(\cdot)$ over
$\Pi_n$ is $+ \infty$.]
To find the minimizer, we use a Lagrange multiplier
$\eta$ and solve the set of equations
\begin{eqnarray*}
\frac{\partial}{\partial\pi_k} g(\pi_1,\ldots, \pi_n;\eta) &=& 0,\qquad 1
\le k \le n,\\
 \frac{\partial}{\partial\eta} g(\pi_1,\ldots, \pi_n;
\eta) &=& 0,
\end{eqnarray*}
where $g(\pi_1,\ldots, \pi_n;\eta)
= \sum_{i=1}^n\log(n \pi_i) - \lambda_n  \sum_{j=1}^{p}\delta_j
( \sum_{i=1}^n\pi_i \Yij )^2 +\break  \eta (\sum_{i=1}^n\pi_i -1)$.
This leads to the equations
\[
{0}= 
\pi^{-1}_k -2 \lambda_n \sum
_{j=1}^{p}\delta_j \Biggl( \sum
_{i=1}^n\pi_i \Yij \Biggr)
\Y_{kj} + \eta, \qquad 1 \le k \le n \quad\mbox{and}\quad {1}= \sum
_{i=1}^n\pi_i,
\]
which, in turn, yield the implicit solution
%
\begin{eqnarray}
\label{pe-empl} \eta&=&2 \lambda_n \sum
_{j=1}^{p}\delta_j M_{nj}^2
-n  \quad \mbox{and}
\nonumber
\\[-8pt]
\\[-8pt]
\nonumber
\pi_k^{-1} &=&  n \Biggl\lbrace1 + 2\gamma_n
\sum_{j=1}^{p}\delta_j
M_{nj}\Y_{kj} - 2 \gamma_n \sum
_{j=1}^{p}M_{nj}^2
\delta_j \Biggr\rbrace, \qquad 1 \le k \le n,
\end{eqnarray}
where $\gamma_n=\lambda_n/n$ and $ M_{nj} = \sum_{i=1}^n\pi_i \Yij$.
To obtain a
more explicit approximation, we show that $\pi_k$'s are of
the form $\pi_k = n^{-1} ( 1 + o_p(1)  ) $
uniformly in $k$. In view of Brouwer's
fixed point theorem [cf. \citet{M65}], it is enough
to show that, with $a_n^{-1} = n^{-1/2}\log(n)$,
%
\begin{eqnarray}
\label{fp-thm} &&\max_{1 \le k \le n} \Biggl\vert\frac{1}{n} \Biggl\lbrace1
 + 2\gamma_n \sum_{j=1}^{p}
\delta_j M_{nj}\Y_{kj} 
 -  2
\gamma_n \sum_{j=1}^{p}
M_{nj}^2\delta_j \Biggr
\rbrace^{-1} - \frac{1}{n} \Biggr\vert
\nonumber
\\[-8pt]
\\[-8pt]
\nonumber
&&\qquad= O_p \bigl(n^{-1}a_n^{-1} \bigr)
\end{eqnarray}
whenever $ \max \lbrace|\pi_k - n^{-1}|
\dvtx1 \le k \le n  \rbrace = O(n^{-1}a_n^{-1})$.
To prove \eqref{fp-thm}, we first show that
%
\begin{equation}
\label{mn-bd} \sup \Biggl\lbrace\gamma_n \sum
_{j=1}^{p}\delta_j \M^2_{nj}
\dvtx  (\pin)' \in\Pi^{0}_n \Biggr
\rbrace = O_p \bigl(a_n^{-2} \bigr),
\end{equation}
where
$\Pi^{0}_n
= \{  ( \pin )' \in\Pi_n \dvtx|\pi_k - n^{-1}| \le
C a_n^{-1}n^{-1}$ {for all} $  1\le k \le n \}$.
Note that for any $ ( \pin )' \in\Pi^{0}_n,  \sum_{i=1}^n\pi_i =1 = \sum_{i=1}^n(1/n)$
$\implies
\sum_{i=1}^n(1 - n \pi_i)=0
\implies\sum_i (1 - n \pi_i)_{+} = \sum_i ( n \pi_i - 1)_{+}
$.
Also, $ (n \pi_i - 1)_{+} > 0$ if and only if (iff)
$ n \pi_i > 1$
and similarly, $ (1 - n \pi_i)_{+} > 0$ { iff} $n \pi_i < 1$.
Hence, using the bound ``$|n \pi_i - 1| \le\C a_n^{-1}$
for all $i=1,\ldots,n$,'' we have
\begin{eqnarray*}
\Biggl\llvert \sum_{i=1}^n(n
\pi_i)^{-1} -n \Biggr\rrvert 
&=& \Biggl\llvert
\sum_{i=1}^n\frac{(1- n \pi_i)_{+}}{n \pi_i} - \sum
_{i=1}^n\frac{( n \pi_i - 1)_{+}}{n \pi_i} \Biggr\rrvert
\\
&=& \Biggl\llvert \sum_{i=1}^n
\frac{(1- n \pi_i)_{+}}{ 1 - (1- n \pi_i)_{+} } - \sum_{i=1}^n
\frac{( n \pi_i - 1)_{+}}{ 1 + ( n \pi_i - 1)_{+}
} \Biggr\rrvert
\\
%
&=& \Biggl\llvert \sum
_{i=1}^n(1- n \pi_i)^2_{+}
\bigl[1 +O ( a_n ) \bigr] + \sum_{i=1}^n(
n \pi_i - 1)^2_{+}   \bigl[1 + O (
a_n ) \bigr] \Biggr\rrvert
\\
&\le& 2 \C^2 n a_n^{-2} \qquad \mbox{for large }
n.
\end{eqnarray*}
%
By  \eqref{pe-empl},
$ n^{-1} \sum_{k=1}^n(n \pi_k)^{-1} - 1=
2 \gamma_n \sum_{j=1}^{p}\delta_j \M_{nj}  [ \bY_{nj} - \M_{nj} ]
$.
Hence by Lem\-ma~\ref{lem-1} and the Cauchy--Schwarz inequality,
\begin{eqnarray*}
\frac{2 \C^2}{a_n^2} &\ge& 2 \gamma_n \Biggl|  \sum
_{j=1}^{p}\delta_j \M_{nj} [
\bY_{nj} - \M_{nj} ]   \Biggr|
\\
&\ge& 2
\Biggl( \gamma_n \sum_{j=1}^{p}
\delta_j \M_{nj}^2 \Biggr)^{1/2}
\Biggl[ \Biggl( \gamma_n \sum_{j=1}^{p}
\delta_j \M_{nj}^2 \Biggr)^{1/2} -
O_p \bigl(n^{-1/2} \bigr) \Biggr]
\end{eqnarray*}
uniformly in $ ( \pin )\in\Pi_n^{0}$, for $n$ large.
Consequently,  \eqref{mn-bd} holds.
Now using  \eqref{mn-bd},  \eqref{pe-empl}, the Cauchy--Schwarz
inequality
and Lemma ~\ref{lem-1},  \eqref{fp-thm} follows.
Hence, by Brouwer's fixed point theorem [cf. \citet{M65}],
there exists
a solution $ ( \pi_1^{0},\ldots,\pi_n^{0} ) $ of  \eqref
{pe-empl} satisfying the bound
%
\begin{equation}
\max \bigl\lbrace\bigl|\pi_k^0 - n^{-1} \bigr| \dvtx1
\le k \le n \bigr\rbrace  = O_p \bigl(n^{-1}a_n^{-1}
\bigr). \label{pi-o-bd}
\end{equation}
Using the second derivative condition, it is easy to verify
that $ ( \pi_1^{0},\ldots,\pi_n^{0} ) $ is a local
minimizer of $f(\cdot)$. In view of the strict
convexity of $f(\cdot)$ on $\Pi_n$, it also follows that $ ( \pi_1^{0},\ldots,\pi_n^{0} ) $ is the
unique minimizer of $f(\cdot)$ over $\Pi_n$.

Next let
$\Ga_{1n} =   2 n \gamma_n^2 \sum_{j=1}^{p}\sum_{l=1}^p
\rho_n(j,l)\delta_j
\delta_l \M_{nj}^{0}\M_{nl}^{0}
$, where
$\M_{nj}^{0}=\break \sum_{i=1}^n\pi_i^{0}\Yij,  1 \le j \le p$.
Then, from  \eqref{pe-empl}, we have
\begin{eqnarray*}
- \log{R}_n(\mu_0)&\equiv& f \bigl(
\pi_1^{0},\ldots,\pi_n^{0} \bigr)
\\
&=& \sum_{i=1}^n\log \Biggl\{ 1 + 2
\gamma_n \sum_{j=1}^{p}
\delta_j \M_{nj}^{0}\Yij - 2
\gamma_n \sum_{j=1}^{p}
\delta_j \bigl(\M_{nj}^{0} \bigr)^2
\Biggr\} \\
&&{}+ n \gamma_n \sum_{j=1}^{p}
\delta_j \bigl(\M_{nj}^{0} \bigr)^2
\\
&\equiv& \sum_{i=1}^n \Biggl\{ 2
\gamma_n \sum_{j=1}^{p}
\delta_j \M_{nj}^{0}\Yij - 2
\gamma_n \sum_{j=1}^{p}
\delta_j \bigl(\M_{nj}^{0} \bigr)^2
\Biggr\} - \Ga_{1n}
\\
&&{} + \R_{1n} + n \gamma_n \sum
_{j=1}^{p}\delta_j \bigl(
\M_{nj}^{0} \bigr)^2
\\
&=& 2 n \gamma_n \sum_{j=1}^{p}
\delta_j \M_{nj}^{0}\bY_{ij}  -  n
\gamma_n \sum_{j=1}^{p}
\delta_j \bigl(\M_{nj}^{0} \bigr)^2
- \Ga_{1n} + \R_{1n}
\\
&\equiv& 2 n \gamma_n \sum_{j=1}^{p}
\M_{nj}^{0}\bY_{ij}  -  n \gamma_n
\sum_{j=1}^{p} \bigl(\M_{nj}^{0}
\bigr)^2
\\
&&{} -2 n \gamma_n^2 \sum_{j=1}^{p}
\sum_{l=1}^p \rho_n(j,l)
\M_{nj}^{0}\M_{nl}^{0} +
\R_{2n}
\\
&\equiv& n \gamma_n ( \bY_{n1}, \ldots,
\bY_{np} ) ( \I_p + 2 \gamma_n \A_n
)^{-1} ( \bY_{n1}, \ldots, \bY_{np}
)'  +  \R_{3n},
\end{eqnarray*}
where the remainder terms $\R_{kn}$'s are defined by subtraction.
By the next lemma, $\R_{kn} = o_p(1)$ for $k=1,2,3$.
Hence, Lemma~\ref{lem-2} is proved.
\end{pf}

\begin{lemma}
\label{lem-3}
Under the conditions of Theorem~\ref{th3.1},
$ \sum_{k=1}^3 |\R_{kn}| = o_p(1)$.
\end{lemma}

\begin{pf}
See [LM] for details.
\end{pf}

\begin{pf*}{Proof of Theorem~\ref{v-str-dep}}
Recall that $\si_{nj}=1,  \forall j,n$. We carry out the proof in
2 steps.\vadjust{\goodbreak}

\textit{Step} (I):  Let
$ \Z_n(t) = \sum_{j=1}^{p}(\sqrt{n} \bY_{nj})
\ind_{ ( {(j-1)}/{p}, {j}/{p}] } (t)$, $ t \in(0,1]$
and let $\Z_n(0) \equiv\Z_n(0+)$. Then,
$\pb (\Z_n \in L^2[0,1]  ) = 1$. The first step is to
prove that $\Z_n(\cdot)\raw^d \Z(\cdot)$ as
elements of $L^2[0,1]$.
Recall that $ \lbrace\ph_j \dvtx j \in\mathbb{N}
\rbrace$ is a complete orthonormal basis for $L^2[0,1]$.
By Theorem 1.84 of Van der Vaart and Wellner
(\citeyear{VW96}), it enough to show that:

(i) For any $0 < t_1 < \cdots< t_r \le1$,  $1 \le r < \infty$,
%
\begin{equation}
\label{NE-1} \bigl( \Z_n(t_1),\ldots,
\Z_n(t_r) \bigr) \rightarrow^{d} \bigl(
\Z(t_1),\ldots, \Z(t_r) \bigr),
\end{equation}

(ii) For any $\varepsilon> 0,  \delta>0 $,
there exists $N=N(\ep,\delta) \in\mathbb{N}$ such that
%
\begin{equation}
\label{NE-2} \limsup_{n \rai}   \pb \Biggl( \sum
_{k=N}^{\infty} \bigl|\langle\Z_n,
\ph_k \rangle\bigr|^2 > \delta \Biggr) < \ep.
\end{equation}
Part (i) can be proved using Theorem 11.1.6 of \citet{AL06}; we omit the routine details.
For part (ii), it is enough to show that
%
\begin{equation}
\label{NET-1} \dm\lim_{N \rai} \limsup_{\nti} \E \Biggl( \sum
_{k=N}^{\infty} \bigl| \langle Z_n,
\phi_k \rangle\bigr|^2 \Biggr) = 0.
\end{equation}
Let $I_j = ( \frac{j-1}{p}, \frac{j}{p} ]$,
$j=1,\ldots,p$. Then, by Fubini's theorem and (C.3),
%
\begin{eqnarray}\label{NET-2}
\E\sum_{k=1}^{\infty} \bigl\llvert \langle
Z_n,\phi_k \rangle \bigr\rrvert^2 &=& \sum
_{k=1}^{\infty}\sum
_{j=1}^{p} \sum_{l=1}^{p}
\int_{I_j} \int_{I_l}
\phi_k(t)\phi_k(s)\rho_0(j/p,l/p) \,\dd s \,
\dd t
\nonumber
\\[-8pt]
\\[-8pt]
\nonumber
&=& \sum_{k=1}^{\infty}\int
_0^1 \int_0^1
\phi_k(t)\phi_k(s)\rho_0(s,t) \,\dd s \,\dd
t +o(1), 
\end{eqnarray}
which equals $\sum_{k=1}^{\infty}\langle\phi_k, \La_0
\phi_k \rangle+o(1)$. Next, using $ |\rho_n(\cdot,\cdot)| \le1$
and $\int|\phi_k(t)| \,\dd t \le
( \int\phi_k^2(t) )^{1/2}=1$, one can show that
for each fixed $k \in\bbn$,
%
\begin{eqnarray} \label{NET-3}
\E \bigl\llvert \langle Z_n,\phi_k \rangle \bigr
\rrvert^2 &=& \sum_{j=1}^{p}
\sum_{l=1}^{p} \int_{I_j}
\int_{I_l} \phi_k(t)\phi_k(s)
\rho_0(j/p,l/p) \,\dd s \,\dd t
\nonumber
\\[-8pt]
\\[-8pt]
\nonumber
&\rightarrow& \int_0^1 \int
_0^1\phi_k(t)\phi_k(s)
\rho_0(s,t) \,\dd s \,\dd t = \langle\phi_k,
\Up_0 \phi_k \rangle.
\end{eqnarray}
By \eqref{NET-2},
\eqref{NET-3} and (C.3), \eqref{NE-2} follows.
Thus, $\Z_n \rightarrow^d \Z$ on $L^2[0,1]$.

%

\textit{Step} (II):  Next we establish weak convergence of the
quadratic form:
%
\begin{equation}
np^{-1} \bbfY_n' ( \I_p + 2
\gamma_n \A_n )^{-1} \bbfY_n = \sum
_{k=0}^{\infty} n p^{-1}
\bbfY_n' (-2 \gamma_n \A_n)^k
\bbfY_n. \label{NE-3}
\end{equation}
Note that by condition (C.3),
$
\Vert\gamma_n \A_n\Vert^2 \le\gamma_n^2 \sum_{j=1}^{p}\sum_{l=1}^{p}
\rho_0^2 ( j/p,l/p  ) \rightarrow\break
c_*^2\int_0^1\int_0^1 \rho_0^2(x,y) \,\dd x \,\dd y    \in(0, 1/4)
$. Hence\vadjust{\goodbreak} there exists a nonrandom $\ep_o \in(0,1)$ such that
for any fixed $m\in\bbn$ and for n sufficiently large
(not depending on~$m$),
%
\begin{eqnarray}\label{NE-4}
&&\sum_{k=m}^{\infty} \bigl\llvert
np^{-1} ( \bY_{n1}, \ldots ,\bY_{np} ) (-2
\gamma_n \A_n)^k ( \bY_{n1},
\ldots,\bY_{np} )' \bigr\rrvert
\nonumber
\\[-8pt]
\\[-8pt]
\nonumber
&&\qquad\le \sum_{k=m}^{\infty} p^{-1}
\sum_{j=1}^{p}(\sqrt{n}
\bY_{nj})^2 \Vert2 \gamma_n \A_n
\Vert^k \le \bigl\Vert\Z_n(\cdot) \bigr\Vert^2  \sum
_{k=m}^{\infty} \eps_o^k.
\end{eqnarray}
Next, let $a_{n.k}(j,l)$ denote the $(j,l)$th element of $(\gamma_n \A_n)^k$.
Note that $\rho_0(\cdot,\cdot)$ is uniformly continuous on $[0,1]^2$.
Hence, using induction and condition (C.3), it can be shown that
for any $k \in\bbn$,
\[
\sum \bigl\{ \bigl| a_{n,k}(j,l) - p^{-1}
c_*^k\rho_0^{*(k)}(j/p,l/p) \bigr|\dvtx 1 \le j,l \le
p \bigr\} = o(1).
\]
By the continuous mapping theorem,
it now follows that for any fixed $m$,
%
\begin{eqnarray}\label{NE-5}
\qquad\hspace*{6pt}\sum_{k=0}^{m} \frac{n}{p}
\bbfY_n'(-2 \gamma_n \A_n)^k
\bbfY_{n} &=& \sum_{k=0}^{m}
\frac{(-2)^k}{p} \sum_{j=1}^{p}\sum
_{k=1}^p \Z_n(j/p)
\Z_n(l/p) a_{n.k}(j,l)
\nonumber
\\[-4pt]
\\[-12pt]
\nonumber
&\rightarrow^d & \sum_{k=0}^{m}
(-2)^k \int_0^1 \int
_0^1 \Z(u) \Z(v) c_*^k
\rho_0^{*(k)}(u,v) \,du \,dv.
\end{eqnarray}
Thus, by \eqref{NE-5}, partial sums of the infinite series
in \eqref{NE-3} converge to the partial sums of the limit series
for any fixed $m$. By \eqref{NE-4}, the tail of the infinite series
in \eqref{NE-3} is negligible. It can be shown (cf. [LM])
that the tail of the limit series is also negligible.
Since, $\Vert\Z_n(\cdot) \Vert_2 \rightarrow^d\Vert\Z(\cdot)
\Vert_2 $,
by \eqref{NE-3}--\eqref{NE-5} and Lemmas~\ref{lem-1},
\ref{lem-3}, the theorem is proved.
\end{pf*}

\subsection{Limit distribution for the ergodic case}\label{sec6.2}
We prove Theorem~\ref{th3.2} and~\ref{th3.3}, using different arguments than
the proof of Theorem~\ref{th3.1}. This is necessitated by the fact that
\textit{we need more accurate bounds on the remainder terms that must
become negligible after the scaling \textup{(}e.g., by $p^{\al_0}$\textup{)}}.
We will also use
the notation $O_p^{u} (\cdot)$ to denote a bound that holds
uniformly over $i \in\{1,\ldots,n\}$ as $\nti$. For example,
$\De_{in} =O_p^{u} (a_n^{-1})$ means
$\max \{ |\De_{in}|\dvtx1\le i\le n  \}= O_p (a_n^{-1})$ an
$\nti$.
Similarly, define $o_p^{u} (\cdot)$.

\begin{pf*}{Proof of Theorem~\ref{th3.2}}
For $1 \le i \le n$, let
$
\De_{in} = 2 \gamma_n \sum_{j=1}^{p}\delta_j \M_{nj}^{0}\Yij$,
$\D_n = 2 \gamma_n \sum_{j=1}^{p}\delta_j (\M_{nj}^{0})^2$ and
$\De_{in}^{0} =  \De_{in} - \D_n$. Then, $\pi_i^0 = 1/[n(1+\De_{in}^0)]$.
Now using $|\pi_i^0 - 1/n| \le a_n^{-1} n^{-1}$ for $1\le i \le n$,
one can show (cf. [LM]) that
%
\begin{equation}
  \bigl|\De_{in}^0\bigr| = O_p^{u}
\bigl(a_n^{-1} \bigr)\quad\mbox{and}\quad |\De_{in}| =
O_p^{u} \bigl(a_n^{-1} \bigr).
\label{del-bd}
\end{equation}
Hence, by Taylor's expansion of $\log(1+x)$ around $x=0$,
%
\begin{eqnarray}\label{EL-E1}
-\log R_n(\mu_0) &=& \sum_{i=1}^n
\log \bigl(1+ \De_{in}^0 \bigr)  +  \lambda\sum
_{j=1}^{p}\delta_j \bigl(M_{nj}^0
\bigr)^2
\nonumber
\\
&=& \Biggl[ 2 n \gamma_n \sum_{j=1}^{p}
\delta_j M_{nj}^0 \bY_{nj} -
\frac
{n}{2}D_n \Biggr] -2^{-1} \Biggl[ \sum
_{i=1}^n\De_{in}^2 +
2D_n \sum_{j=1}^n
\De_{jn} \Biggr]
\\
&&{} + 3^{-1}\sum_{i=1}^n
\De_{in}^3  + O_p \bigl(na_n^{-4}
\bigr).\nonumber
\end{eqnarray}
There exist
$E_{nj} = {n^{-1}\sum_{i=1}^n}[Y_{ij} \cdot O^u_p(a_n^{-4})]$, $1\leq
j\leq p$,
such that (cf. [LM])
%
\begin{eqnarray}\label{M-exp}
M_{nj}^0 
&=&
\bY_{nj} - {n^{-1}\sum_{i=1}^n}Y_{ij}
\De_{in} + D_n\bY_{nj} + {n^{-1}\sum
_{i=1}^n}Y_{ij} \bigl[
\De_{in}^2 +2 \De_{in}D_n \bigr]
\nonumber
\\[-8pt]
\\[-8pt]
\nonumber
&&{} -{n^{-1}\sum_{i=1}^n}Y_{ij}
(\De_{in})^3 + E_{nj}; 
\end{eqnarray}
\begin{eqnarray*}
L_{1n}&\equiv& 2 \gamma_nn \sum
_{j=1}^{p} \delta_jM_{nj}^0
\bY_{nj} - \frac{n}{2}D_n \qquad\mbox{the lead term of}
- \log R_n(\mu)
\\
%
&=& \gamma_nn \sum
_{j=1}^{p}\delta_jM_{nj}^0
\bY_{nj} + \sum_{i=1}^n
\De_{in}^2/2 + \sum_{i=1}^n
\De_{in}^4/2 -\frac{1}{2} \Biggl( \sum
_{i=1}^n\De_{in} \Biggr)D_n
\\
&&{} - \frac{1}{2}\sum_{i=1}^n
\De_{in}^3 - \sum_{i=1}^n
\De_{in}^2D_n + \sum
_{i=1}^n\De_{in} \cdot
O^u_p \bigl(a_n^{-4} \bigr).
\end{eqnarray*}
%
Hence, from \eqref{EL-E1},
%
\begin{eqnarray}\label{EL-E2}
-\log R_n(\mu) &=& n \gamma_n \sum
_{j=1}^{p}\delta_jM_{nj}^0
\bY_{nj} -\frac{3}{2} \Biggl( \sum_{i=1}^n
\De_{in} \Biggr)D_n - \frac{1}{6} \sum
_{i=1}^n\De_{in}^3
\nonumber
\\[-8pt]
\\[-8pt]
\nonumber
&&{}+ O_p \bigl(n a_n^{-4} \bigr).
\end{eqnarray}
Next, define $\be_n^2 = p \sum_{k=1}^p k^{-\al}$. Then,
%
\begin{equation}
\beta_n^2\sim\cases{ %
 Cp,
&\quad $\mbox{if } \al>1,$
\vspace*{2pt}\cr
Cp\log p, &\quad $\mbox{if } \al=1,$
\vspace*{2pt}\cr
Cp^{2-\al}, & \quad $\mbox{if } 0< \al<1. $}
\label{betan}
\end{equation}
Note that for any $1 \le i \le n$,
$\E\bY_{nj} Y_{ij}=n^{-1}\sum_{k=1}^n \E Y_{kj}Y_{ij} = 1/n$. Let
$V_{in}
\equiv
\gamma_n\sum_{j=1}^{p}(\delta_j\bY_{nj} Y_{ij}-1/n)$
and $D_{nj}^Y \equiv\{y\dvtx P(Y_{1j} = y) >0\}$.
Then there exists a $\delta\in(0,1)$
such that by using \eqref{delta-ex} and the conditions of
Theorem~\ref{th3.2}, one gets (cf. [LM])
%
\begin{eqnarray}\label{Ibd-E1}
\E\sum_{i=1}^nV_{in}^2
&\le& C \Biggl[n \E \Biggl\{\gamma_n\sum
_{j=1}^{p}(\bY_{nj} Y_{1j} - 1/n)
\Biggr\}^2 +n \E \Biggl\{ \gamma_n\sum
_{j=1}^{p}\delta_j \bY_{nj}^3
Y_{1j} \Biggr\}^2
\nonumber
\\[-2pt]
&&{} + n \E \Biggl\{ \gamma_n \sum_{j=1}^{p}
\delta_j\bW_n(j,j)\bY_{nj} Y_{1j}
\Biggr\}^2 \nonumber\\[-2pt]
&&{}+ n\E \Biggl\{ \gamma_n \sum
_{j=1}^{p} Y_{ij}^2
\ind(s_{nj}=0) \Biggr\}^2 \Biggr]
\\[-2pt]
&\leq& C p^{-2}\be_n^2 + O
\bigl(n^{-1} \bigr) + n \delta^{n-2} \sum
_{y \in D^Y_{nj}} y^2 \pb(Y_{1j}=y)
\nonumber
\\[-2pt]
&\leq& C p^{-2} \be_n^2 + O
\bigl(n^{-1} \bigr).\nonumber
\end{eqnarray}
Using \eqref{M-exp}, \eqref{Ibd-E1} and the Cauchy--Schwartz inequality,
one can show that
%
\begin{eqnarray}\label{EL-E3}
n \gamma_n \sum_{j=1}^{p}
\delta_jM_{nj}^0\bY_{nj} &=& n
\gamma_n \sum_{j=1}^{p}
\delta_j\bY_{nj}^2 - \sum
_{i=1}^n\De_{in}[V_{in} + 1/n]
\nonumber
\\[-9pt]
\\[-9pt]
\nonumber
&&{} +O_p \bigl(a_n^{-2} \bigr) +
O_p \bigl(a_n^{-1}p^{-1}
\be_n \bigr).
\end{eqnarray}
Also, from \eqref{M-exp}, for any $1\le k \le n$, we have
\begin{eqnarray*}
\De_{kn} &=& 2 \gamma_n \sum
_{j=1}^{p}\delta_jM_{nj}^0Y_{kj}
\\[-2pt]
&\equiv& 2 \gamma_n \sum_{j=1}^{p}
\delta_j \bY_{nj} Y_{kj} - \frac{2}{n}
\De_{kn} + R_{1n}(k)
\\[-2pt]
&\equiv& \biggl( 1 + \frac{2}{n} \biggr)^{-1} 2
\gamma_n \Biggl( \sum_{j=1}^{p}
\bY_{nj} Y_{kj} \Biggr) + R_{2n}(k),
\end{eqnarray*}
where $R_{1n}(k)$ and $R_{2n}(k)$ are remainder terms
satisfying (cf. [LM])
\begin{eqnarray*}
\sum_{k=1}^n R_{ln}^2(k)
&=& O_p \bigl( n a_n^{-2} p^{-2}
\be_n^2 + n a_n^{-4} \bigr),\qquad
l=1,2. \label{Ibd-E3}
\end{eqnarray*}
Next, using similar arguments
and noting that $\E\bW_n(j,j) \bY_{nj}^2 =O(n^{-2})$ and
$\Var(\bW_n(j,j) \bY_{nj}^2) \le C n^{-3}$ for all $j$,
one can show (cf. [LM]) that
%
\begin{eqnarray}\qquad
\sum_{i=1}^n\De_{in}
(V_{in} + 1/n) &= &\biggl(\frac{n}{n+2} \biggr) 2
\gamma_n^2 \sum_{i=1}^n
\Biggl( \sum_{j=1}^{p}\bY_{nj}
Y_{ij} \Biggr)^2 \bigl(1+ o_p(1) \bigr),
\label{Ibd-E4}
\\[-2pt]
 n \gamma_n \sum_{j=1}^{p}
\delta_j \bY_{nj}^2 &=& n \gamma_n
\sum_{j=1}^{p}\bY_{nj}^2
+ O_p \bigl(n^{-1} + n^{-1/2} p^{-1}
\be_n \bigr), \label{Ibd-E5}
\\[-2pt]
\label{Ibd-E6}\sum_{i=1}^n\De_{in}^3
&= &\biggl(\frac{n}{n+2} \biggr) 8 \gamma_n^2 \Biggl[
\sum_{i=1}^n \Biggl( \sum
_{j=1}^{p}\bY_{nj} Y_{ij}
\Biggr)^3 \Biggr] \bigl(1+o_p(1) \bigr)
\nonumber
\\[-8pt]
\\[-8pt]
\nonumber
&=& O_p \bigl(n^{-1} + n^{-1/2} p^{-1}
\be_n \bigr) .
\end{eqnarray}
Using \eqref{EL-E2}, \eqref{EL-E3} and \eqref{Ibd-E4}--\eqref{Ibd-E6}
and the fact that $\E\bY_{nj} Y_{ij}= n^{-1}$ and $\Var(\bY_{nj}
Y_{ij} )
\le C n^{-1}$ for all $i,j$, one can conclude (cf. [LM])
%
\begin{equation}
-\log R_n(\mu_0) = n \gamma_n \sum
_{j=1}^{p}\bY_{nj}^2 -
\frac{2}{n} + O_p \bigl(n a_n^{-4} \bigr)
+ O_p \bigl(a_n^{-1}p^{-1}
\be_n \bigr). \label{EL-E4}
\end{equation}
Set $a_n=n^{1/2}(p/n^2)^{1/10}$. Then $a_n \rai$, $a_n=o(n^{1/2})$
and $\sqrt{p}.n. a_n^{-4}=o(1),$ there by making the last two terms
in \eqref{EL-E4} $o(p^{-1/2})$ whenever $p= o(n^2)$. Now Theorem~\ref{th3.2}
follows by adapting the proof of the CLT for a stationary sequence
of $\varrho$-mixing random variables to triangular arrays.
\end{pf*}

\begin{pf*}{Proof of Theorem~\ref{th3.3}} Arguments
in the proof of Theorem~\ref{th3.2} yield
the asymptotic approximation for $-\log R_n(\mu_0)$ in \eqref{EL-E4}
with $\be_n^2 \sim c p^{2-\al}$ for all $0 < \al< 1/2$.
Now choose $a_n \rai$ to satisfy
$
p^{\al} [{n}^{-1} + n a_n^{-4} + a_n^{-1}p^{-1}\be_n]
 \rightarrow  0
$,
for example, $a_n=n^{1/2}  [{p^{\al}}/{n} ]^{1/5}$.
Then it follows that
\[
p^{\al} \bigl( -\log R_n(\mu_0) - c_* \bigr) =
c_* p^{\al} \Biggl( \frac{1}{p} \sum
_{j=1}^{p} \bigl\{ (\sqrt {n}\bY_{nj})^2
-1 \bigr\} \Biggr) + o_p(1).
\]
In the case where $X_{nj}$'s are Gaussian with $\varrho_n(i,j)
=\varrho_{\al}(i-j)$, the leading term has the same distribution as
$W_n \equiv (p^{-2+2\al} )^{1/2} \sum_{j=1}^{p}(Z_j^2 - 1)$,
where $\{Z_j\}$ is a stationary Gaussian process with correlation
function $\varrho_{\al}(\cdot)$. Then the result of Taqqu (\citeyear{T75})
implies the $W_n \raw^d W$, and the theorem follows. In the general
case when $X_{nj}$'s are not Gaussian,
the theorem follows by using
convergence of moments of $\sqrt{n}\bY_n$ to the
moments of $N(0,1)$ and
a variant of the diagram formula; cf. \citet{A94}.
\end{pf*}

\begin{pf*}{Proof of Theorem~\ref{th3.4}}
Similar to the proof of Theorem~\ref{th3.3}; see [LM].
\end{pf*}

\begin{pf*}{Proof of Theorem~\ref{th4.1}} Note that
$ \cI_n= \{I_i \dvtx1 \le i \le n-m+1 \}$ with
$I_i=\{i,i+1,\ldots, i+m-1\}$. Let
\[
U_j(x) \equiv\dm\sum_{i=(j-1)m+1}^{jm \land n -m+1}
\ind \bigl( -\log R_m^{*}(\mu_0,I_i)
\le x \bigr),  \qquad  x\in\bbr,
\]
for $1\leq j\leq M$
where $M= \lceil(n-m+1)/m \rceil$ is the smallest integer not
less than $(n-m+1)/m$. By the independence of $U_j(x)$
and $U_{j+k}(x)$ for $k \ge2$, one can show (cf. [LM]) that
for each
$x \in\bbr$,
%
\begin{eqnarray}\label{sub-cg}
&&\E \bigl( \hG_{n}^{\mathrm{NE}}(x) -\pb \bigl( -\log
R_m^{*}(\mu_0,I_i) \le x \bigr)
\bigr)^2
\nonumber
\\[-8pt]
\\[-8pt]
\nonumber
&&\qquad\le C n^{-2} M m^2 = o(1).
\end{eqnarray}
The next arguments are similar to the proof of the
Glivenko--Cantelli
theorem [cf. Theorem 13.3 of Billingsley (\citeyear{Bill99})]
and the continuity of
the limit distribution of $-\log R_n(\mu_0)$; one can complete
the proof; see [LM].
\end{pf*}

\begin{pf*}{Proof of Theorem~\ref{th4.2}}  Similar to the
proof of Theorem~\ref{th4.1}.
\end{pf*}

\begin{pf*}{Proof of Remark 4.1}  Here we outline a
proof of consistency of the permutation invariant estimators
$\hat{\al}$ and $\hat{\ka}$ of Remark 4.1. W.l.g, suppose that
$\mu=0$ and $\si_{nj}=0$ for all $j$, $n$.
First consider\vspace*{1pt} the estimator $\hat{\al}$ in \eqref{p-inv-hal}.
Let $A_n =\{s_{nj}\not=0$ for all $j=1,\ldots,p\}$, $n\geq1$.
Then, by condition (C.1), $P(A_n^c) = O(p \delta^{n-1})$ for some
$\delta\in(0,1)$. On the set $A_n$, $e_n=0$, and
using arguments similar to those in the proof of Lemma~\ref{lem-2},
one can show that
\begin{eqnarray*}
e_n+ n^{-1}\sum_{i=1}^n
\Biggl\{p^{-1}\sum_{j=1}^p
[X_{ij} - \bX_{nj}]\delta_j^{1/2}
\Biggr\}^2 = n^{-1}\sum_{i=1}^n
\bX_{ip}^2 + R_n,
\end{eqnarray*}
where $R_n = o_p(p^{-\al})$.
Note that $EI_{11} \sim\frac{2Cp^{-\al}}{(1-\al)(2-\al)}$
for $0<\al<1$ and $\operatorname{Var}(I_{11})=O(n^{-1}p^{-2\al})$. Now
it is easy to verify that $\hal$ satisfies \eqref{al-cond}.

To prove the consistency of $\hat{\ka}$ of
\eqref{p-inv-hal}, using moderate deviation inequalities
[cf. G\"otze and Hipp (\citeyear{GH78})], one
can
conclude that
%
\begin{equation}\qquad\label{mod-inq}
P \bigl(\bigl|\hat{c}(j,k) - c(j,k)\bigr| > n^{-1/2}\log n \mbox{ for some } 1
\leq j,k\leq p \bigr)
= o(1).
\end{equation}
Next note that $Ck^{-\al}\leq n^{-1/2}\log n$
for all $k > n^{{1}/{(2\al)}}$. This implies that
only $O(n^{{1}/{(2\al)}})$-many correlation terms contribute
to $\hat{\ka}$, with probability tending to one. Now using~\eqref{mod-inq}
for the nonvanishing terms and the fact that $\al>1/2$,
one can prove consistency of $\hat{\ka}$.\vspace*{-2pt}
\end{pf*}

\section*{Acknowledgments}\vspace*{-5pt}
We thank three anonymous reviewers and
Professors
Peter B\"uhlmann and Art Owen
for a number of constructive comments that improved
an earlier version of the paper.\vspace*{-2pt}

\begin{supplement}
\stitle{Numerical results and proofs}
\slink[doi]{10.1214/12-AOS1040SUPP} 
\sdatatype{.pdf}
\sfilename{aos1040\_supp.pdf}
\sdescription{Additional simulation results and some details of proofs.}
\end{supplement}

%

\printaddresses

\end{document}